\documentclass[11pt]{article}
\usepackage{amsmath,amssymb}

\newcommand{\epsln}{\varepsilon}

\def\epsln{\varepsilon}
\def\bgeq{\begin{equation}}
\def\edeq{\end{equation}}
\def\bgar{\begin{array}}
\def\edar{\end{array}}

\title{Asymptotic Stability of the Stationary Solution for a Hyperbolic
  Free Boundary Problem Modeling Tumor Growth}
\author{Shangbin Cui}
\date{\small Institute of Mathematics, Sun Yat-Sen University, Guangzhou,
  Guangdong 510275,\\ People's Republic of China. E-mail:\,cuisb3@yahoo.com.cn}

\begin{document}

\maketitle

\begin{abstract}
  In this paper we study asymptotic behavior of solutions for a free boundary
  problem modeling the growth of tumors containing two species of cells:
  proliferating cells and quiescent cells. This tumor model was proposed by
  Pettet et al in {\em Bull. Math. Biol.} (2001). By using a functional
  approach and the $C_0$ semigroup theory, we prove that the unique stationary
  solution of this model ensured by the work of Cui and Friedman ({\em Trans.
  Amer. Math. Soc.}, 2003) is locally asymptotically stable in certain function
  spaces. Key techniques used in the proof include an improvement of the
  linear estimate obtained by the work of Chen et al ({\em Trans. Amer. Math.
  Soc.}, 2005), and a similarity transformation.

\medskip

   {\bf Keywords and phrases}: Free boundary problem; hyperbolic equations;
   tumor growth; stationary solution; asymptotic stability.
\medskip

   {\bf AMS subject classification:} 35C10, 35Q80, 92C15.

\end{abstract}

\section{Introduction}
\setcounter{equation}{0}

  During the past thirty years, an increasing number of free boundary problems
  of partial differential equations have been proposed by groups of researchers
  to model the growth of various {\em in vivo} and {\em in vitro} tumors, see,
  e.g., \cite{tumrev}, \cite{ByrCha1}--\cite{ByrCha3}, \cite{Franketal},
  \cite{Fried}, \cite{Green1}, \cite{Green2}, \cite{Pettet}--\cite{WarKin} and
  the references cited therein. Such free boundary problems usually contain one
  or more reaction diffusion equations describing the distribution of nutrient
  and inhibitory materials, and several first-order nonlinear partial
  differential equations or nonlinear conservation laws with source terms
  describing the evolution and movement of various tumor cells (proliferating
  cells, quiescent cells and dead cells). Rigorous analysis of such tumor
  models is evidently a significant topic of research and has drawn great
  attention during the past a few years. Main concern of this topic is the
  dynamics or the long-term behavior of solutions of such free boundary
  problems.

  Based on applications of the well-established theories of elliptic and
  parabolic partial differential equations, parabolic differential equations
  in Banach spaces (i.e., differential equations in Banach spaces that are
  treatable with the analytic semigroup theory) and the bifurcation theory,
  rigorous analysis of models for the growth of tumors containing only one
  species of tumor cells has achieved great depth, cf. \cite{BazFri1},
  \cite{BazFri2}, \cite{Cui1}, \cite{Cui2}, \cite{Cui4}--\cite{CuiEsc},
  \cite{FriHu1}--\cite{FriRei2}, \cite{WuCui}, \cite{ZhouCui} and the
  references cited therein. As far as models for tumors containing more than
  one species of tumor cells are concerned, however, the progress is relatively
  backward. This is caused by the fact that such tumor models are much more
  difficult to analyze because they contain nonlinear conservation laws whose
  dynamical behavior is very hard to grasp.

  In this paper we study the following free boundary problem modeling the
  growth of an {\em in vitro} tumor containing two species of cells ---
  proliferating cells and quiescent cells:
\begin{equation}
   \nabla^2 C=F(C) \quad \mbox{for} \;\; x\in\Omega(t),\;\; t\geq 0,
\end{equation}
%---(1.1)---
\begin{equation}
   C=C_0\quad \mbox{for} \;\; x\in\partial\Omega(t),\;\; t\geq 0,
\end{equation}
%---(1.2)---
\begin{equation}
  {\partial P\over\partial t}+\nabla\cdot(\vec{u}P)=\big[K_B(C)-K_Q(C)\big]P
  +K_P(C)Q \quad \mbox{for} \;\; x\in\Omega(t),\;\; t\geq 0,
\end{equation}
%---(1.3)---
\begin{equation}
  {\partial Q\over\partial t}+\nabla\cdot(\vec{u}Q)=K_Q(C)P-\big[K_D(C)+K_P(C)
  \big]Q \quad \mbox{for} \;\; x\in\Omega(t),\;\; t\geq 0,
\end{equation}
%---(1.4)---
\begin{equation}
    P+Q=N \quad \mbox{for} \;\; x\in\Omega(t),\;\; t\geq 0,
\end{equation}
%---(1.5)---
\begin{equation}
   {dR\over dt}=\vec{u}\cdot\vec{\nu}\quad \mbox{for} \;\;
   x\in\partial\Omega(t),\;\; t\geq 0.
\end{equation}
%---(1.6)---
  Here $C$ denotes the concentration of nutrient (with all nutrient materials
  regarded as one species), $P$ and $Q$ denote the densities of proliferating
  cells and quiescent cells, respectively, whose mixture makes up the tumor
  tissue and has a constant density $N$, $\vec{u}$ denotes the velocity of the
  cell movement, $R$ denotes the radius of the tumor, $\Omega(t)=\{x\in
  {\mathbb R}^3: r=|x|<R(t)\}$ is the domain occupied by the tumor at time $t$,
  and $\vec{\nu}$ is the unit outward normal of $\partial\Omega(t)$.
  Besides, $C_0$ is a positive constant reflecting the constant nutrient supply
  that the tumor receives from its surface, $F(C)$ is the nutrient
  consumption rate function, and $K_B(C)$, $K_D(C)$, $K_P(C)$ and $K_Q(C)$ are
  the birth rate of proliferating cells, death rate of quiescent cells,
  transferring rate of proliferating cells to quiescent cells and transferring
  rate of quiescent cells to proliferating cells, respectively. We shall
  only consider radially symmetric solutions of the above problem, so that $C$,
  $P$, $Q$ are functions of the radial space variable $r=|x|$ and the time
  variable $t$, and $\vec{u}=u(r,t)r^{-1}x$, where $u$ is a scaler function.

  The above tumor model was proposed by Pettet {\em et al} in the literature
  \cite{Pettet}. Its global well-posedness has been established by Cui and
  Friedman in \cite{CuiFri1}. A challenging task concerning this free boundary
  problem is the study of the asymptotic behavior of its solutions as time
  goes to infinity. For the corresponding model of the growth of tumors with
  one species of cells, it is known that there exists a unique stationary
  solution and all time-dependent solutions converge to it as time goes to
  infinity, or in other words, this unique stationary solution is globally
  asymptotically stable, cf. \cite{Cui1} and \cite{FriRei1}. Since the above
  problem is a natural extension of such one species tumor model to the two
  species case, we are naturally lead to the conjecture that a similar result
  holds for it. Advancement of the study toward this goal is as follows. In
  \cite{CuiFri2}, Cui and Friedman proved that the problem (1.1)--(1.6) has a
  unique stationary solution. In \cite{ChenCuiF}, Chen, Cui and Friedman
  further proved that this stationary solution is linearly asymptotically
  stable, namely, the trivial solution of the linearization of (1.1)--(1.6) at
  the stationary solution is asymptotically stable. However, this
  last-mentioned result does not imply, at least straightforwardly, that the
  stationary solution of (1.1)--(1.6) is asymptotically stable. In fact, to
  the best of our knowledge this problem has been remaining open before this
  manuscript is prepared. We refer the reader to see \cite{ChenFri},
  \cite{Cui2} and \cite{CuiWei} for other related work.

  In this paper we shall prove that the unique stationary solution of
  (1.1)--(1.6) ensured by  \cite{CuiFri2} is locally asymptotically stable.
  Recall that conditions given in \cite{CuiFri2} which ensure that (1.1)--(1.6)
  has a unique stationary solution are as follows:
\begin{equation}
  F(C),\;\;K_B(C),\;\;K_D(C),\;\;K_P(C)\;\;{\rm and}\;\;K_Q(C)\;\;{\rm are
  \;\; analytic\;\; in}\;\;C,\;\;0\leq C\leq C_0;
\end{equation}
%---(1.7)---
\begin{equation}
  F(0)=0,\quad F'(C)>0 \quad {\rm for}\;\;0\leq C\leq C_0;
\end{equation}
%---(1.8)---
\begin{equation}
\left\{
\begin{array}{l}
   K_B'(C)>0\;\;\mbox{and}\;\;K_D'(C)<0\;\;{\rm for}\;\;0\leq C\leq C_0,\;\;
   K_B(0)=0\;\;{\rm and}\;\;K_D(C_0)=0;\\
   K_P(C)\;\;\mbox{and}\;\;K_Q(C)\;\;\mbox{satisfy the same conditions as}
   \;\; K_B(C)\;\;\mbox{and}\;\;K_D(C),\;\;\mbox{respectively};\\
   K_B'(C)+K_D'(C)>0\;\;{\rm for}\;\; 0\leq C\leq C_0.
\end{array}
\right.
\end{equation}
%---(1.9)---
  The main result of this paper is the following:
\medskip

  {\bf Theorem 1.1}\ \ {\em Assume that the conditions $(1.7)$--$(1.9)$
  are satisfied. Let $(C_*,P_*,Q_*,\vec{u}_*,R_*)$ be the unique stationary
  solution of the problem $(1.1)$--$(1.6)$, and let $(C,P,Q,\vec{u},R)$ be a
  time-dependent solution of it such that $P|_{t=0}=P_0$, $Q|_{t=0}=Q_0$ and
  $R|_{t=0}=R_0$, where $P_0$, $Q_0$ and $R_0$ are given initial data
  satisfying $0\leq P_0\leq N$, $0\leq Q_0\leq N$ and $P_0+Q_0=N$. Then there
  exist positive constants $\mu$, $\epsln$ and $K$ such that if $P_0$, $Q_0$
  and $R_0$ satisfy
$$
  \max_{0\leq r\leq 1}|P_0(rR_0)-P_*(rR_*)|<\epsln, \quad
  \sup_{0<r<1}r(1\!-\!r)\Big|{dP_0(rR_0)\over dr}-
  {dP_*(rR_*)\over dr}\Big|<\epsln,
$$
$$
  \max_{0\leq r\leq 1}|Q_0(rR_0)-Q_*(rR_*)|<\epsln, \quad
  \sup_{0<r<1}r(1\!-\!r)\Big|{dQ_0(rR_0)\over dr}-
  {dQ_*(rR_*)\over dr}\Big|<\epsln
$$
  and $|R_0-R_*|<\epsln$, then for all $t\geq 0$ we have
$$
  \max_{0\leq r\leq 1}|P(rR(t),t)-P_*(rR_*)|<K\epsln e^{-\mu t}, \quad
  \sup_{0<r<1}r(1\!-\!r)\Big|{\partial P(rR(t),t)\over\partial r}-
  {dP_*(rR_*)\over dr}\Big|<K\epsln e^{-\mu t},
$$
$$
  \max_{0\leq r\leq 1}|Q(rR(t),t)-Q_*(rR_*)|<K\epsln e^{-\mu t}, \quad
  \sup_{0<r<1}r(1\!-\!r)\Big|{\partial Q(rR(t),t)\over\partial r}-
  {dQ_*(rR_*)\over dr}\Big|<K\epsln e^{-\mu t}
$$
  and $|R(t)-R_*|<K\epsln e^{-\mu t}$.}
\medskip

  We shall use a functional approach to prove the above theorem. More precisely,
  we shall first reduce the problem (1.1)--(1.6) into a differential equation
  for the unknown $U=(p,z)$ in the Banach space $X=C[0,1]\times {\mathbb R}$,
  where $p=p(r,t)=P(rR(t),t)$ and $z=z(t)=\log R(t)$. The reduced equation is
  of the hyperbolic type in the sense of Pazy \cite{Pazy}, and is quasi-linear.
  We next use the Banach fixed point theorem to prove that for any $U_0=(p_0,
  z_0)$ sufficiently closed to the stationary point $U_*=(p_*,z_*)$, where $p_*
  =p_*(r)=P_*(rR_*)$ and $z_*=\log R_*$, this differential equation imposed
  with the initial condition $U|_{t=0}=U_0$ and the decay estimate $\sup_{t\geq
  0}e^{\mu t}\|U(t)-U_*\|_{X_0}<\infty$, where $X_0$ is a subspace of $X$, has
  a unique solution in the space $C([0,\infty),X_0)$ (endowed with the norm
  $|\|U|\|=\sup_{t\geq 0}e^{\mu t}\|U(t)\|_{X_0}$). To attain this goal we shall
  use some abstract results for hyperbolic differential equations in Banach
  spaces established in \cite{Pazy}. In particular, a family of evolution
  systems for the linear equations related to the semi-linearization of the
  reduced equation are obtained and applied to convert the semi-linearized
  equations into integral equations. The main difficult and key step in the
  proof of Theorem 1.1 is the establishment of a uniform decay estimate for the
  family of evolution systems. To obtain it, we first use a localization
  technique to get an improvement of the linear estimate established in
  \cite{ChenCuiF}, removing the singularities at $r=0$ contained in that
  estimate. See Lemma 6.2 in Section 6 for this improved linear estimate. We
  next develop a similarity transformation technique and use it to extend this
  improved linear estimate to the family of evolution systems mentioned above.
  See Section 5 for details of this transformation and Lemma 6.4 in Section 6
  for the uniform decay estimate for the family of evolution systems.

  The layout of the rest part is as follows. In the following section we reduce
  the problem (1.1)--(1.6) into a differential equation in the Banach space
  $X=C[0,1]\times {\mathbb R}$. In Section 3 we summarize some basic properties
  of the stationary solution. The reader is suggested to pay attention to
  properties of the stationary solution at the end point $r=0$ which will play
  an important role in later analysis. In Section 4 we prove that the linear
  parts of semi-linearizations of the reduced equation are related with a stable
  family of generators of $C_0$ semigroups on $X$, so that their solution
  operators are evolution systems. This result enables us to use the abstract
  results of \cite{Pazy} to convert the semi-linearized equations into integral
  equations. The most important technique used in this paper --- similarity
  transformations --- will be developed in Section 5. In Section 6 we first
  derive an improvement of the linear estimate established in \cite{ChenCuiF}
  and next use the similarity transformation technique to extend this estimate
  to the evolution systems obtained in Section 4. After these preparations, in
  the last section we use the Banach fixed point theorem to prove Theorem 1.1.

  Throughout this paper the notation ``$\,'\,$'' denotes both the ordinary
  derivatives of functions in ${\mathbb R}$ and the Fr\'{e}chet derivatives of
  mappings between Banach spaces.

\section{Reduction of the problem}
\setcounter{equation}{0}

  In this section we reduce the system of equations (1.1)--(1.6) into a
  differential equation in the Banach space $X=C[0,1]\times {\mathbb R}$.

  We first note that by summing up (1.3), (1.4) and using (1.5), we get the
  following equation:
\begin{equation}
    \nabla\cdot\vec{u}={1\over N}\big[K_B(C)P-K_D(C)Q\big].
\end{equation}
%---(2.1)---
  Conversely, from (1.3), (1.5) and (2.1) we immediately obtain (1.4). Hence,
  the two groups of equations (1.3), (1.4), (1.5) and (1.3), (1.5), (2.1) are
  equivalent.

  By rescaling the space and time variables, setting
$$
   p={P\over N},\quad q={Q\over N}=1-p,\quad c={C\over C_0},\quad
   \vec{u}=u{x\over|x|},
$$
  and using the equivalence of (1.3), (1.4) and (1.5) with (1.3), (1.5) and
  (2.1), we see that the problem (1.1)--(1.6) can be reformulated into the
  following form:
\begin{equation}
   {\partial^2 c\over\partial r^2}+{2\over r}{\partial c\over\partial r}
   =F(c) \quad \mbox{for} \;\; 0<r\leq R(t),\;\; t\geq 0,
\end{equation}
%---(2.2)---
\begin{equation}
   {\partial c\over\partial r}\Big|_{r=0}=0,\quad c|_{r=R(t)}=1 \quad
   \mbox{for}\;\; t\geq 0,
\end{equation}
%---(2.3)---
\begin{equation}
   {\partial p\over\partial t}+u{\partial p\over\partial r}=
   K_P(c)\!+\!\big[K_M(c)\!-\!K_N(c)\big]p\!-\!K_M(c)p^2 \quad
   \mbox{for} \;\; 0\leq r\leq R(t),\;\; t\geq 0,
\end{equation}
%---(2.4)---
\begin{equation}
   {\partial u\over\partial r}+{2\over r}u=-K_D(c)+K_M(c)p \quad
   \mbox{for} \;\; 0<r\leq R(t)\;\; \mbox{and}\;\;u|_{r=0}=0,\;\; t\geq 0,
\end{equation}
%---(2.5)---
\begin{equation}
   {dR\over dt}=u(R,t)  \quad \mbox{for} \;\; t\geq 0,
\end{equation}
%---(2.6)---
  where
\begin{equation}
   K_M(c)=K_B(c)+K_D(c),\quad K_N(c)=K_P(c)+K_Q(c),
\end{equation}
%---(2.7)---
  and $F(c)$, $K_B(c)$, $K_D(c)$, $K_P(c)$ and $K_Q(c)$ are rescaled forms of
  the corresponding functions appearing in (1.1)--(1.6).

  Next, we set
$$
   \bar{c}(\bar{r},t)=c(\bar{r}e^{z(t)},t), \quad \bar{p}(\bar{r},t)=
   p(\bar{r}e^{z(t)},t), \quad  \bar{u}(\bar{r},t)=
   u(\bar{r}e^{z(t)},t)e^{-z(t)}, \quad R(t)=e^{z(t)},
$$
  where $0\leq\bar{r}\leq 1$, $t\geq 0$. Then the problem (2.2)--(2.6) is
  further reduced into the following problem (for simplicity of the notation
  we omit all bar's):
\begin{equation}
   {\partial^2 c\over\partial r^2}+{2\over r}{\partial c\over\partial r}
   =e^{2z} F(c) \quad \mbox{for} \;\; 0<r\leq 1,\;\; t\geq 0,
\end{equation}
%---(2.8)---
\begin{equation}
   {\partial c\over\partial r}\Big|_{r=0}=0,\quad c|_{r=1}=1 \quad
   \mbox{for}\;\; t\geq 0,
\end{equation}
%---(2.9)---
\begin{eqnarray}
   \displaystyle{\partial p\over\partial t}+[u(r,t)-ru(1,t)]
   {\partial p\over\partial r}=K_P(c)+&&\big[K_M(c)\!-\!K_N(c)\big]p
   -\!K_M(c)p^2 \nonumber\\ [0.2cm]
   && \mbox{for} \;\; 0\leq r\leq 1,\;\; t\geq 0,
\end{eqnarray}
%---(2.10)---
\begin{equation}
   {\partial u\over\partial r}+{2\over r}u=-K_D(c)+K_M(c)p \quad
   \mbox{for} \;\; 0<r\leq 1\;\; \mbox{and}\;\;u|_{r=0}=0,\;\; t\geq 0,
\end{equation}
%---(2.11)---
\begin{equation}
   {dz\over dt}=u(1,t)  \quad \mbox{for} \;\; t\geq 0.
\end{equation}
%---(2.12)---

  To further reduce (2.8)--(2.12) we first note that (2.8) and (2.9) can be
  solved to express $c$ as a function of $z$. Thus, instead of $c(r,t)$, later
  on we shall use the notation $c(r,z(t))$ or simply $c(r,z)$ to denote the
  solution of (2.8) and (2.9). Next, we note that (2.11) can be solved to get
  $u$ as a functional of $p$ and $z$. Thus later on we use the notation
  $u_{p,z}$ to re-denote $u$. By a simple computation we have
\begin{equation}
  u_{p,z}(r,t)={1\over r^2}\int_0^r [-K_D(c(\rho,z(t)))+
   K_M(c(\rho,z(t)))p(\rho,t)]\rho^2 d\rho
\end{equation}
%---(2.13)---
  for $0<r\leq 1$, $t\geq 0$, and $u_{p,z}(0,t)=0$ for $t\geq 0$. We also
  denote
\begin{equation}
  w_{p,z}(r,t)=u_{p,z}(r,t)-ru_{p,z}(1,t).
\end{equation}
%---(2.14)---
  It follows that (2.8)--(2.12) reduces into the following system of equations :
\begin{equation}
\left\{
\begin{array}{l}
   \displaystyle{\partial p\over\partial t}+w_{p,z}(r,t)
   {\partial p\over\partial r}=f(r,p,z)\quad
   \mbox{for} \;\; 0\leq r\leq 1,\;\; t>0,\\
   \displaystyle{dz\over dt}=u_{p,z}(1,t)  \quad \mbox{for} \;\; t>0,
\end{array}
\right.
\end{equation}
%---(2.15)---
  where
$$
  f(r,p,z)=K_P(c(r,z))+\big[K_M(c(r,z))\!-\!K_N(c(r,z))\big]p-K_M(c(r,z))p^2.
$$
\medskip

  In what follows we shall rewrite (2.15) as a differential equation in the
  Banach space $X=C[0,1]\times {\mathbb R}$. Let
$$
  C^1_V[0,1]=\{p\in C[0,1]\cap C^1(0,1):\; r(1-r)p'(r)\in C[0,1]\},
$$
  with norm
$$
  \|p\|_{C^1_V[0,1]}=\max_{0\leq r\leq 1}|p(r)|+
  \sup_{0<r<1}|r(1-r)p'(r)|\quad \mbox{for} \;\; p\in C^1_V[0,1].
$$
  It is evident that $C^1_V[0,1]$ endowed with this norm is a Banach space
  densely and continuously embedded into $C[0,1]$. Given $p\in C[0,1]$ and
  $z\in {\mathbb R}$, we introduce a linear operator ${\mathcal A}_0(p,z):
  C^1_V[0,1]\to C(0,1)$ as follows: For any $q\in C^1_V[0,1]$,
$$
  {\mathcal A}_0(p,z)q(r)=-w_{p,z}(r)q'(r) \quad \mbox{for} \;\; 0<r<1.
$$
  Here and hereafter $w_{p,z}(r)$ represents the function defined by similar
  formulations as in (2.13) and (2.14), with $p(r,t)$ and $z(t)$ there replaced
  by $p(r)$ and $z$, respectively.
  Later on we shall use the convention that for a function $f\in C(0,1)$, if
  both limits $\lim_{r\to 0^+}f(r)$ and $\lim_{r\to 1^-}f(r)$ exist and are
  finite, then we write $f\in C[0,1]$. Furthermore, when we are concerned with
  the values of $f$ at $r=0$ and $r=1$, we mean that $f(0)=\lim_{r\to 0^+}f(r)$
  and $f(1)=\lim_{r\to 1^-}f(r)$. Using this convention, we see easily that for
  any $p\in C[0,1]$ and $z\in {\mathbb R}$ we have $w_{p,z}(r)/r(1\!-\!r)\in
  C[0,1]$. It follows that for any $q\in C^1_V[0,1]$,
  both limits $\lim_{r\to 0^+}w_{p,z}(r)q'(r)$ and $\lim_{r\to 1^-}w_{p,z}(r)
  q'(r)$ exist, so that ${\mathcal A}_0(p,z)q\in C[0,1]$. It can also be easily
  seen that ${\mathcal A}_0(p,z)$ is a bounded linear operator from
  $C^1_V[0,1]$ to $C[0,1]$, and
$$
  \|{\mathcal A}_0(p,z)\|_{L(C^1_V[0,1],C[0,1])}\leq \sup_{0<r<1}
  \Big|{w_{p,z}(r)\over r(1\!-\!r)}\Big|.
$$
  Next we introduce mappings ${\mathcal F}:C[0,1]\times {\mathbb R}\to C[0,1]$
  and ${\mathcal G}:C[0,1]\times {\mathbb R}\to {\mathbb R}$ respectively by
$$
  {\mathcal F}(p,z)(r)=f(r,p(r),c(r,z)),
$$
  and
$$
  {\mathcal G}(p,z)=\displaystyle\int_0^1 [-K_D(c(r,z))+K_M(c(r,z))p(r)]r^2 dr.
$$
  We set $X_0=C^1_V[0,1]\times {\mathbb R}$, which is a Banach
  space with the product norm and is densely and continuously embedded into $X=
  C[0,1]\times {\mathbb R}$. We now define a nonlinear operator ${\mathbb F}:
  X_0\to X$ as follows:
$$
  {\mathbb F}(U)=\big({\mathcal A}_0(p,z)p+{\mathcal F}(p,z),
  {\mathcal G}(p,z)\big)  \quad \mbox{for}\;\;U=(p,z)\in X_0.
$$
  It is obvious that ${\mathbb F}\in C^\infty(X_0,X)$. Later on we shall also
  regard ${\mathbb F}$ as an unbounded nonlinear operator in $X$ with domain
  $X_0$. With these notation and convention, we can rewrite (2.15) as the
  following differential equation in the Banach space $X$:
\begin{equation}
  {d U\over dt}={\mathbb F}(U).
\end{equation}
%---(2.16)---
  Here $U=U(t)$ represents a $X_0$-valued unknown function for $t\geq 0$, and
  the left-hand side denotes the Fr\'{e}chet derivative of $U=U(t)$ regarded as
  a mapping from $[0,\infty)$ to the $X$ space.
\medskip

  It will be convenient to denote, for $U=(p,z)\in X$ and $V=(q,y)\in X_0$,
$$
   {\mathbb A}_0(U)V=\big({\mathcal A}_0(p,z)q, 0\big) \quad \mbox{and} \quad
  {\mathbb F}_0(U)=\big({\mathcal F}(p,z), {\mathcal G}(p,z)\big).
$$
  Then we have
$$
  {\mathbb F}(U)={\mathbb A}_0(U)U+{\mathbb F}_0(U)
  \quad \mbox{for}\;\; U\in X_0.
$$
  Clearly, for every $U\in X$, ${\mathbb A}_0(U)$ is a bounded linear operator
  from $X_0$ to $X$, i.e, ${\mathbb A}_0(U)\in L(X_0,X)$. Furthermore, it
  can be easily seen that ${\mathbb A}_0\in C^\infty(X,L(X_0,X))$. Later on we
  shall also regard ${\mathbb A}_0(U)$ as an unbounded linear operator in $X$
  with domain $X_0$. Finally, we note that ${\mathbb F}_0\in C^\infty(X,X)\cap
  C^\infty(X_0,X_0)$.
\medskip

  From \cite{CuiFri2} we know that under the conditions (1.7)--(1.9) which we
  assume to be true throughout the whole paper, the problem (2.2)--(2.6) has a
  unique stationary solution. It follows that the problem (2.8)--(2.12) has a
  unique stationary solution which we denote as $(c_*,p_*,u_*,z_*)$. By
  definition, $(c_*,p_*,u_*,z_*)=(c_*(r),p_*(r),u_*(r),z_*)$ ($0\leq r\leq 1$)
  is the solution of the following problem:
\setcounter{equation}{16}
\begin{equation}
   c_*''+{2\over r}c_*'=e^{2z_*} F(c_*) \quad \mbox{for} \;\; 0<r\leq 1,
\end{equation}
%---(2.17)---
\begin{equation}
   c_*'(0)=0,\quad c_*(1)=1,
\end{equation}
%---(2.18)---
\begin{eqnarray}
   u_* p_*'=f(r,p_*,z_*) \quad \mbox{for} \;\; 0\leq r\leq 1,
\end{eqnarray}
%---(2.19)---
\begin{equation}
   u_*'+{2\over r}u_*=-K_D(c_*)+K_M(c_*)p_* \quad
   \mbox{for} \;\; 0<r\leq 1,
\end{equation}
%---(2.20)---
\begin{equation}
   u_*(0)=0, \quad  u_*(1)=0.
\end{equation}
%---(2.21)---
  Let $U_*=(p_*,z_*)$. Then $U_*\in X_0$ (see Lemma 2.1 below) and it is the
  unique equilibrium of (2.16), i.e.,
$$
  {\mathbb F}(U_*)=0,
$$
  or
$$
  {\mathbb A}_0(U_*)U_*+{\mathbb F}_0(U_*)=0.
$$

  Since our goal is to study asymptotic stability of the stationary solution
  $U_*$, it will be convenient to rewrite (2.16) into an equation for the
  difference $V=U-U_*$. For this purpose we introduce two nonlinear operators
  ${\mathbb A}:X\to L(X_0,X)$ and ${\mathbb G}: X\to X$ as follows:
$$
  {\mathbb A}(V)W={\mathbb A}_0(U_*+V)W+[{\mathbb A}_0'(U_*)W]U_*
  +{\mathbb F}_0'(U_*)W \quad \mbox{for}\;\; V\in X,\;\; W\in X_0,
$$
$$
\begin{array}{rl}
  {\mathbb G}(V)=&[{\mathbb A}_0(U_*+V)-{\mathbb A}_0(U_*)-
  {\mathbb A}_0'(U_*)V]U_* \\
  &+[{\mathbb F}_0(U_*+V)-{\mathbb F}_0(U_*)
  -{\mathbb F}_0'(U_*)V] \quad \mbox{for}\;\; V\in X.
\end{array}
$$
  Then clearly (2.16) can be rewritten as the following equivalent equation
  for $V=U-U_*$:
$$
  {d V\over dt}={\mathbb A}(V)V+{\mathbb G}(V),
\eqno{(2.22)}
$$
  i.e., if $U$ is a solution of (2.16) then $V=U-U_*$ is a solution of (2.22)
  and vice versa. We note that ${\mathbb A}\in C^\infty(X,L(X_0,X))$,
  ${\mathbb G}\in C^\infty(X,X)$, and by using the Taylor expansions up to
  second-order for Fr\'{e}chet derivatives of ${\mathbb A}_0$ and ${\mathbb
  F}_0$ we have
$$
  \|{\mathbb G}(V)\|_X=O(\|V\|_X^2) \quad \mbox{as}\;\; \|V\|_X\to 0.
\eqno{(2.23)}
$$
  We also note that, by introducing an operator ${\mathbb B}: X\to X$ by
$$
  {\mathbb B}W=[{\mathbb A}_0'(U_*)W]U_*
  +{\mathbb F}_0'(U_*)W \quad \mbox{for}\;\; W\in X,
$$
  we have
$$
  {\mathbb A}(0)={\mathbb F}'(U_*)={\mathbb A}_0(U_*)+{\mathbb B}
  \quad \mbox{and} \quad
  {\mathbb A}(V)={\mathbb A}_0(U_*+V)+{\mathbb B}.
$$
  Note that as an immediate consequence of the facts that ${\mathbb A}\in
  C^\infty(X,L(X_0,X))$ and ${\mathbb F}_0\in C^\infty(X,X)$, we have
  ${\mathbb B}\in L(X)$. We also note that $[V\to {\mathbb A}_0(U_*+V)]\in
  C^\infty(X,L(X_0,X))$.

  From the above deduction it follows immediately that the stationary solution
  $(c_*,p_*,u_*,z_*)$ of (2.8)--(2.12) is asymptotically stable if and only if
  the trivial solution of (2.22) is asymptotically stable. More precisely,
  Theorem 1.1 follows if we prove that the solution $V=V(t)$ of (2.22)
  satisfies $\|V(t)\|_{X_0}\leq K\epsln e^{-\mu t}$, $t\geq 0$, provided
  $\|V(0)\|_{X_0}\leq\epsln$ for some small $\epsln>0$. Hence, later on we
  shall concentrate our attention on the equation (2.22).

  A simple computation shows that if we denote
$$
  a(r)=K_M(c_*(r))-K_N(c_*(r))-2K_M(c_*(r))p_*(r),
\eqno{(2.24)}
$$
\vspace*{-0.2cm}
\setcounter{equation}{24}
\begin{eqnarray}
  b(r)&=&\{K_P'(c_*(r))+[K_M'(c_*(r))-K_N'(c_*(r))]p_*(r)-K_M'(c_*(r))p_*^2(r)\}
  c_z(r) \nonumber\\[0.2cm]
  &&\displaystyle+rp_*'(r)\Big[\int_0^1 g_c(\rho)c_z(\rho)\rho^2 d\rho-
  {1\over r^3}\int_0^r g_c(\rho)c_z(\rho)\rho^2 d\rho\Big],
\end{eqnarray}
%---(2.25)---
$$
  {\mathcal B}q(r)=rp_*'(r)\Big[\int_0^1 g_p(\rho)q(\rho)\rho^2 d\rho-
  {1\over r^3}\int_0^r g_p(\rho)q(\rho)\rho^2 d\rho\Big],
\eqno{(2.26)}
$$
$$
  {\mathcal F}(q)=\int_0^1 g_p(\rho)q(\rho)\rho^2 d\rho,
\eqno{(2.27)}
$$
  and
$$
  \kappa=\int_0^1 g_c(\rho)c_z(\rho)\rho^2 d\rho,
\eqno{(2.28)}
$$
  where
$$
  g_p(r)=K_M(c_*(r)), \quad g_c(r)=-K_D'(c_*(r))+K_M'(c_*(r))p_*(r),
  \quad c_z(r)={\partial c\over\partial z}(r,z_*),
$$
  then we have
$$
  {\mathbb B}=
\left(
\begin{array}{cc}
  a(r)+{\mathcal B}\;\; & \;\; b(r) \\ {\mathcal F} \;\; & \;\; \kappa
\end{array}
\right).
\eqno{(2.29)}
$$
  Here and hereafter, when we write ${\mathbb M}=\displaystyle\left(
  \begin{array}{cc} M_{11}\;\; & \;\; M_{12} \\ M_{21} \;\; & \;\; M_{22}
  \end{array}\right)$ for bounded linear operators $M_{11}\in L(X_1,Y_1)$,
  $M_{12}\in L(X_2,Y_1)$, $M_{21}\in L(X_1,Y_2)$, $M_{22}\in L(X_2,Y_2)$,
  where $X_1$, $X_2$, $Y_1$ and $Y_2$ are Banach spaces, we mean that
  ${\mathbb M}$ is the bounded linear operator from $X_1\times X_2$ to
  $Y_1\times Y_2$ defined by
$$
  {\mathbb M}(x_1,x_2)=(M_{11}x_1+M_{12}x_2,M_{21}x_1+M_{22}x_2) \quad
  \mbox{for}\;\; (x_1,x_2)\in X_1\times X_2.
$$
  Using this notation we see that
$$
    {\mathbb A}_0(U_*)=
\left(
\begin{array}{cc}
  {\mathcal L}_0 \;\; & \;\; 0 \\ 0 \;\; & \;\; 0
\end{array}
\right), \quad \mbox{where}\;\; {\mathcal L}_0={\mathcal A}_0(p_*,z_*),
$$
  and, for $V=(\varphi,\zeta)\in X$,
$$
  {\mathbb A}_0(U_*+V)=
\left(
\begin{array}{cc}
  {\mathcal L}_V\;\; & \;\; 0 \\ 0 \;\; & \;\; 0
\end{array}
\right), \quad \mbox{where}\;\; {\mathcal L}_V=
  {\mathcal A}_0(p_*+\varphi,z_*+\zeta).
\eqno{(2.30)}
$$
  We recall that $a(r)<0$ for all $0\leq r\leq 1$, see (2.7) in Section 2 of
  \cite{ChenCuiF}.

\section{Some basic facts}
\setcounter{equation}{0}

  We summarize some basic properties of the functions $c_*(r)=c(r,z_*)$,
  $c_z(r)=\displaystyle{\partial c\over\partial z}(r,z_*)$, $p_*(r)$ and
  $u_*(r)$ in the following lemma. These properties will play an important
  role in later discussions.
\medskip

  {\bf Lemma 3.1}\ \ {\em We have the following assertions:

  $(1)$ $c_*,\, c_z\in C^\infty[0,1]$, and
$$
  0<c_*(0)\leq c_*(r)\leq 1 \quad \mbox{for}\;\; 0\leq r\leq 1, \quad
  c_*'(r)>0 \quad \mbox{for}\;\; 0<r\leq 1, \quad c_*'(0)=0.
\eqno{(3.1)}
$$

  $(2)$ $p_*\in C[0,1]\cap C^\infty(0,1]$,
$$
  0<p_*(0)\leq p_*(r)\leq 1 \quad \mbox{for}\;\; 0\leq r\leq 1, \quad
  p_*'(r)>0\;\;\;\mbox{for}\;\; 0<r\leq 1,
\eqno{(3.2)}
$$
  and either $p_*\in C^1[0,1]$ or there exists $0<\gamma<1$ such that
  $\lim_{r\to 0^+}r^\gamma p_*'(r)$ exists and is finite, so that $r^\gamma
  p_*'(r)\in C[0,1]$. Moreover,
$$
  \lim_{r\to 0^+}rp_*'(r)=0, \quad \lim_{r\to 0^+}r^2p_*''(r)=0, \quad
  \lim_{r\to 0^+}r^3p_*'''(r)=0.
\eqno{(3.3)}
$$

  $(3)$ $u_*\in C^1[0,1]\cap C^\infty(0,1]$, and there exist positive constants
  $C_1$, $C_2$ such that
$$
  -C_1r(1-r)\leq u_*(r)\leq -C_2r(1-r) \quad \mbox{for}\;\;\; 0\leq r\leq 1.
\eqno{(3.4)}
$$
  Besides, either $u_*\in C^2[0,1]$ or there exists $0<\gamma<1$ such that
  $\lim_{r\to 0^+}r^\gamma u_*''(r)$ exists and is finite, so that
  $r^\gamma u_*''(r)\in C[0,1]$. Moreover,
$$
  \lim_{r\to 0^+}ru_*''(r)=0, \quad \lim_{r\to 0^+}r^2u_*'''(r)=0.
\eqno{(3.5)}
$$
}
\medskip

  {\em Proof}:\ \ The assertions that $c_*,c_z\in C^\infty[0,1]$ and relations
  in (3.1) are immediate. The assertions that $p_*\in C[0,1]\cap C^\infty(0,1]$,
  $u_*\in C^1[0,1]\cap C^\infty(0,1]$ and relations in (3.2) follow from
  Theorem 2.1 of \cite{CuiFri2}, by which we also know that $u_*(r)<0$ for
  $0<r<1$. The last assertion combined with the facts that $u_*'(0)<0$ and
  $u_*'(1)>0$ (see Theorem 7.1 of \cite{CuiFri2}) immediately yields (3.4). To
  prove (3.3) we compute:
$$
  \lim_{r\to 0}u_*(r)p_*'(r)=K_P(c_*(0))+[K_M(c_*(0))-K_N(c_*(0))]p_*(0)
  -K_M(c_*(0))p_*^2(0)=0
$$
  (see (8.4) in Section 8 of \cite{CuiFri2}), so that
$$
  \lim_{r\to 0}rp_*'(r)=\lim_{r\to 0}{r\over u_*(r)}\cdot
  \lim_{r\to 0}u_*(r)p_*'(r)=0,
$$
  and
$$
\begin{array}{rcl}
  \displaystyle\lim_{r\to 0}u_*(r)rp_*''(r)&=&\displaystyle\lim_{r\to 0}r
  \{K_P'(c_*(r))+[K_M'(c_*(r))-K_N'(c_*(r))]p_*(r)
  \\[0.2cm]
  &&\displaystyle-K_M'(c_*(r))p_*^2(r)\}c_*'(r)+\lim_{r\to 0}
  \{[K_M(c_*(r))-K_N(c_*(r))]
  \\[0.2cm]
  &&\displaystyle-2K_M(c_*(r))p_*(r)\}rp_*'(r)-\lim_{r\to 0}
  u_*'(r)rp_*'(r)=0
\end{array}
$$
  so that
$$
  \lim_{r\to 0}r^2p_*''(r)=\lim_{r\to 0}{r\over u_*(r)}\cdot
  \lim_{r\to 0}u_*(r)rp_*''(r)=0.
$$
  This proves the first two relations in (3.3). The proof of the third relation
  is similar and is omitted. Next, from (2.20) we can easily deduce that
$$
\begin{array}{rcl}
  \displaystyle\qquad\quad
  u_*''(r)&=&[-K_D'(c_*(r))+K_M'(c_*(r))p_*(r)]c_*'(r)
  +K_M(c_*(r))p_*'(r) \\ [0.2cm]
  && +\displaystyle{2\over r}[K_D(c_*(r))-K_M(c_*(r))p_*(r)]+
  {6\over r^4}\int_0^r[-K_D(c_*(\rho))+K_M(c_*(\rho))p_*(\rho)]\rho^2 d\rho
  \\  [0.2cm]
  &=&\displaystyle[-K_D'(c_*(r))+K_M'(c_*(r))p_*(r)]c_*'(r)+
  {6\over r^4}\int_0^r[K_D(c_*(r))-K_D(c_*(\rho))]\rho^2 d\rho \\ [0.2cm]
  &&+K_M(c_*(r))p_*'(r)-\displaystyle
  {6\over r^4}\int_0^r[K_M(c_*(r))p_*(r)-K_M(c_*(\rho))p_*(\rho)]\rho^2 d\rho.
  \hfill (3.6)
\end{array}
$$
  From this expression and the first relation in (3.3) we readily obtain
  the first relation in (3.5). The proof of the second relation in (3.5) is
  similar and is omitted. Finally, by Theorems 5.3 and 5.4 of \cite{CuiFri2}
  we know that either $p_*\in C^1[0,1]$ or there exist constants $-1<\alpha<0$
  and $C$ such that\footnotemark[1]
\footnotetext[1]{In the notation of Theorem 5.4 of \cite{CuiFri2}, we have
  $\alpha=\alpha(\lambda)$ and $C=(1+\alpha(\lambda))\omega$.}
$$
  p_*'(r)=Cr^\alpha+O(1) \quad \mbox{for}\;\; r\to 0.
\eqno{(3.7)}
$$
  Suppose that it is the second case. Then, by letting $\gamma=|\alpha|$, we
  see that $0<\gamma<1$ and $r^\gamma p_*'(r)\in C[0,1]$. Finally, from (3.7)
  we see that
$$
  p_*(r)=p_*(0)+C(1+\alpha)^{-1}r^{1+\alpha}+O(r) \quad \mbox{for}\;\; r\to 0.
\eqno{(3.8)}
$$
  Substituting (3.7) and (3.8) into (3.6) we get $r^\gamma u_*''(r)\in C[0,1]$.
  This completes the proof of Lemma 3.1. $\quad\Box$
\medskip

  {\bf Corollary 3.2}\ \ {\em Let $a$, $b$, ${\mathcal B}$, ${\mathcal F}$ and
  ${\mathbb B}$ be as in $(2.24)$--$(2.27)$ and $(2.29)$. Then we have $a,b\in
  C^1_V[0,1]$, ${\mathcal B}\in L(C[0,1],C^1_V[0,1])\subseteq L(C[0,1])\cap
  L(C^1_V[0,1])$, ${\mathcal F}\in L(C[0,1],{\mathbb R})$, and ${\mathbb B}\in
  L(X)\cap L(X_0)$. Moreover, we also have $r^2(1\!-\!r)^2 a''(r),\,
  r^2(1\!-\!r)^2 b''(r)\in C[0,1]$. $\quad\Box$}
\medskip

  {\bf Corollary 3.3}\ \ {\em ${\mathbb G}\in C^\infty(X,X)\cap C^\infty(X_0,
  X_0)$, and in addition to $(2.23)$ we also have:
$$
  \|{\mathbb G}(V)\|_{X_0}=O(\|V\|_{X_0}^2) \quad \mbox{as}\;\;
  \|V\|_{X_0}\to 0.
\eqno{(3.9)}
$$
}

  {\em Proof}:\ \ We have ${\mathbb G}(V)={\mathbb G}_1(V)+{\mathbb G}_2(V)$,
  where
$$
  {\mathbb G}_1(V)=[{\mathbb A}_0(U_*+V)-{\mathbb A}_0(U_*)-
  {\mathbb A}_0'(U_*)V]U_*,
$$
$$
  {\mathbb G}_2(V)={\mathbb F}_0(U_*+V)-{\mathbb F}_0(U_*)
  -{\mathbb F}_0'(U_*)V.
$$
  Since ${\mathbb F}_0\in C^\infty(X_0,X_0)$, it is evident that ${\mathbb G}_2
  \in C^\infty(X_0,X_0)$ and $\|{\mathbb G}_2(V)\|_{X_0}=O(\|V\|_{X_0}^2)$ as
  $\|V\|_{X_0}\to 0$. Next, let $V=(\varphi,\zeta)$ and $(p,z)=(p_*+\varphi,z_*
  +\zeta)$. Then by (2.30) we have
$$
  {\mathbb A}_0(U_*+V)U_*=(-w_{p,z}(r)p_*'(r),0).
$$
  Using this expression and the first two relations in (3.3) we can easily show
  that for every $V\in X$ we have ${\mathbb A}_0(U_*+V)U_*\in X_0$, and the
  mapping $V\to {\mathbb A}_0(U_*+V)U_*$ belongs to $C^\infty(X,X_0)$. Hence we
  have ${\mathbb G}_1\in C^\infty(X,X_0)\subseteq C^\infty(X_0,X_0)$ and
  $\|{\mathbb G}_1(V)\|_{X_0}=O(\|V\|_{X}^2)=O(\|V\|_{X_0}^2)$ as $\|V\|_{X_0}
  \to 0$. Combining these assertions together, we see that the desired
  assertion follows.
  $\quad\Box$

\section{Evolution systems}
\setcounter{equation}{0}

  Given a small positive number $\epsln$, we denote
$$
  S_\epsln=\{V=(\varphi,\zeta)\in X=C[0,1]\times {\mathbb R}:\;
  \|\varphi\|_\infty\leq\epsln,\; |\zeta|\leq\epsln\}.
$$
  In this section we shall prove that the family of operators $\{{\mathbb A}
  (V):V\in S_\epsln\}$ is a stable family of infinitesimal generators of $C_0$
  semigroups on $X$, and its part in $X_0$ is a stable family of infinitesimal
  generators of $C_0$ semigroups on $X_0$. For the concept of {\em stable
  family of infinitesimal generators of $C_0$ semigroups} and related results,
  we refer the reader to see Sections 5.2--5.5, Chapter 5 and Section 6.4,
  Chapter 6 of \cite{Pazy}. We use the notation $\tilde{\mathbb A}(V)$ to
  denote the part of ${\mathbb A}(V)$ in $X_0$. Recall that
$$
  {\rm Dom}(\tilde{\mathbb A}(V))=\{U\in X_0:{\mathbb A}(V)U\in X_0\},
  \quad \mbox{and}
$$
$$
  \tilde{\mathbb A}(V)U={\mathbb A}(V)U \quad \mbox{for}\;\;
  U\in {\rm Dom}(\tilde{\mathbb A}(V)).
$$

  Let $w\in C^1[0,1]$ and assume that it satisfies the following condition:
  There exist positive constants $C_1$ and $C_2$ such that
$$
  -C_1r(1-r)\leq w(r)\leq -C_2r(1-r) \quad \mbox{for}\;\; 0\leq r\leq 1.
\eqno{(4.1)}
$$
  Note that this assumption particularly implies that $w(0)=w(1)=0$, $w'(0)<0$
  and $w'(1)>0$. For a such $w\in C^1[0,1]$, we denote by ${\mathcal L}_0$ the
  bounded linear operator from $C^1_V[0,1]$ to $C[0,1]$ defined by
$$
  {\mathcal L}_0q(r)=-w(r)q'(r) \quad \mbox{for} \;\; 0<r<1, \quad
  \mbox{for}\;\; q\in C^1_V[0,1].
$$
  Later on we shall also regard ${\mathcal L}_0$ as an unbounded linear
  operator in $C[0,1]$ with domain $C^1_V[0,1]$. Note that if $w=u_*$ then
  ${\mathcal L}_0={\mathcal A}_0(p_*,z_*)$.
\medskip

  {\bf Lemma 4.1}\ \ {\em Let the notation and the assumption be as above. Then
  ${\mathcal L}_0$ generates a $C_0$ semigroup of contractions $e^{t{\mathcal
  L}_0}$ on $C[0,1]$, i.e.,
$$
  \|e^{t{\mathcal L}_0}\|_{L(C[0,1])}\leq 1 \quad \mbox{for}\;\; t\geq 0.
\eqno{(4.2)}
$$
  Moreover, $C^1_V[0,1]$ is ${\mathcal L}_0$-admissible\footnotemark[2],
\footnotetext[2]{Recall that for a $C_0$ semigroup $T(t)$ ($t\geq 0$) on a
  Banach space $X$ generated by an unbounded linear operator $A$ in $X$, a
  linear subspace $Y$ of $X$ is called {\em $A$-admissible} if it is an
  invariant subspace of $T(t)$  for all $t\geq 0$, and the restriction of
  $T(t)$ ($t\geq 0$) to $Y$ is a $C_0$ semigroup in $Y$. A necessary and
  sufficient condition for $Y$ to be $A$-admissible is that (1) $Y$ is an
  invariant subspace of $R(\lambda, A)$ for all $\lambda>\omega$ and (2) the
  part $\bar{A}$ of $A$ in $Y$ is an infinitesimal generator of a $C_0$
  semigroup on $Y$. In this case we have $e^{t\bar{A}}=e^{tA}|_Y$. See
  Theorem 5.5 in Chapter 4 of \cite{Pazy}.}
  and the restriction of $e^{t{\mathcal L}_0}$ on $C^1_V[0,1]$ is a uniformly
  bounded $C_0$ semigroup on $C^1_V[0,1]$, i.e., there exists constant $C>0$
  depending only on the constants $C_1$ and $C_2$ in $(4.1)$ such that
$$
  \|e^{t{\mathcal L}_0}\|_{L(C^1_V[0,1])}\leq C \quad \mbox{for}\;\; t\geq 0.
\eqno{(4.3)}
$$
}

  {\em Proof}:\ \ We first prove that for any $\lambda\in {\mathbb C}$ with
  ${\rm Re}\lambda>0$ and any $f\in C[0,1]$, the equation
$$
  -w(r)q'(r)-\lambda q(r)=f(r)
\eqno{(4.4)}
$$
  has a unique solution $q\in C^1_V[0,1]$, and $\|q\|_\infty\leq ({\rm Re}
  \lambda)^{-1}\|f\|_\infty$.

  Arbitrarily take a number $0<r_0<1$ and fix it. Since the equation (4.4)
  is linear and regular for all $0<r<1$, for each given $c\in {\mathbb R}$ it
  has a unique solution for all $0<r<1$ satisfying $q(r_0)=c$. In fact, this
  solution is given by
$$
  q(r)=e^{-\lambda\int_{r_0}^r {d\rho\over w(\rho)}}\Big[c-
  \int_{r_0}^r {f(\eta)\over w(\eta)} e^{\lambda\int_{r_0}^{\eta}
  {d\rho\over w(\rho)}}d\eta\Big].
\eqno{(4.5)}
$$
  Since $w\in C^1[0,1]$ and $w(0)=0$, we have $w(r)=w'(0)r[1+o(1)]=w'(0)r
  [1+o(1)]^{-1}$ for $r\sim 0^+$. Thus, by taking $\delta>0$ sufficiently
  small, we see that
$$
  \int_{r_0}^r{d\rho\over w(\rho)}={1\over w'(0)}\int_\delta^r
  {1+o(1)\over\rho}d\rho+\int_{r_0}^\delta{d\rho\over w(\rho)}
  ={1+o(1)\over w'(0)}\log r+C
\eqno{(4.6)}
$$
  for $r\sim 0^+$. Since $w'(0)<0$ and ${\rm Re}\lambda>0$, it follows that
$$
  e^{-\lambda\int_{r_0}^r {d\rho\over w(\rho)}}=
  Cr^{\lambda [1+o(1)]\over |w'(0)|}\to 0 \quad \mbox{as}\;\; r\to 0^+.
$$
  Hence, using the L'Hospital's law we see that for any $c\in {\mathbb R}$
  the function $q(r)$ given by (4.5) has finite limit as $r\to 0^+$. By a
  similar argument as in the deduction of (4.6) we have
$$
  e^{-\lambda\int_{r_0}^r {d\rho\over w(\rho)}}=
  C(1-r)^{-{\lambda [1+o(1)]\over w'(1)}}\to \infty
   \quad \mbox{as}\;\; r\to 1^-.
$$
  It follows that the function $q(r)$ given by (4.5) cannot be bounded in a
  neighborhood of $r=1$ unless we take $c=\displaystyle\int_{r_0}^1 {f(\eta)
  \over w(\eta)}e^{\lambda\int_{r_0}^{\eta}{d\rho\over w(\rho)}}
  d\eta$, which gives
$$
  q(r)=e^{-\lambda\int_{r_0}^r {d\rho\over w(\rho)}}
  \int_r^1 {f(\eta)\over w(\eta)} e^{\lambda\int_{r_0}^{\eta}
  {d\rho\over w(\rho)}} d\eta.
\eqno{(4.7)}
$$
  By using the L'Hospital's law we can easily verify that the function $q(r)$
  given by (4.7) has finite limit as $r\to 1^-$. Hence, we have shown that
  if ${\rm Re}\lambda>0$ then for any $f\in C[0,1]$ the equation (4.4) has a
  unique solution $q\in C[0,1]\cap C^1(0,1)$. Using (4.4) as well as the fact
  that $w(r)/r(1\!-\!r)\in C[0,1]$ we see readily that $q\in C^1_V[0,1]$.
  Furthermore, by a simple computation we have
$$
\begin{array}{rcl}
  |q(r)|&\leq &\displaystyle\|f\|_\infty
  e^{({\rm Re}\lambda)\int_{r_0}^r {d\rho\over |w(\rho)|}}\int_r^1
  {1\over |w(\eta)|}
  e^{-({\rm Re}\lambda)\int_{r_0}^\eta {d\rho\over |w(\rho)|}}d\eta
\\ [0.3cm]
  &=&\displaystyle {\|f\|_\infty\over {\rm Re}\lambda}
  \Big[1-e^{-({\rm Re}\lambda)\int_r^1 {d\rho\over |w(\rho)|}}\Big]
  \leq {\|f\|_\infty\over {\rm Re}\lambda} \quad
  \mbox{for}\;\; 0<r<1.
\end{array}
$$
  This proves the desired assertion.

  Thus, we have proved that for any $\lambda\in {\mathbb C}$ with ${\rm Re}
  \lambda>0$, there holds $\lambda\in\rho({\mathcal L}_0)$ and
$$
  \|R(\lambda,{\mathcal L}_0)\|_{L(C[0,1])}\leq {1\over {\rm Re}\lambda}.
\eqno{(4.8)}
$$
  It follows by the Hille-Yosida Theorem that ${\mathcal L}_0$ generates a
  strongly continuous semigroup $e^{t{\mathcal L}_0}$ on $C[0,1]$ which
  satisfies the estimate (4.2).

  Next we prove that $C^1_V[0,1]$ is ${\mathcal L}_0$-admissible. Clearly,
  $C^1_V[0,1]$ is an invariant subspace of $R(\lambda,{\mathcal L}_0)$ for all
  $\lambda\in {\mathbb C}$ with ${\rm Re}\lambda>0$, because we know that for
  such $\lambda$, $R(\lambda,{\mathcal L}_0)$ is a bounded linear operator in
  $C[0,1]$ with image contained in ${\rm Dom}({\mathcal L}_0)=C^1_V[0,1]$. Let
  $\tilde{\mathcal L}_0$ be the part of ${\mathcal L}_0$ in $C^1_V[0,1]$. Since
  $r(1\!-\!r)/w(r),\, w(r)/r(1\!-\!r)\in C[0,1]$, we see that for $q\in C[0,1]
  \cap C^1(0,1)$, $r(1\!-\!r)q'(r)\in C[0,1]$ if and only if $w(r)q'(r)\in
  C[0,1]$, and, furthermore, there exist positive constants $C_1$ and $C_2$
  such that
$$
  C_1\sup_{0<r<1}|w(r)q'(r)|\leq \sup_{0<r<1}|r(1\!-\!r)q'(r)|\leq
  C_2\sup_{0<r<1}|w(r)q'(r)|.
$$
  By the above assertion it follows easily that
$$
  {\rm Dom}(\tilde{\mathcal L}_0)=
  \{q\in C[0,1]\cap C^2(0,1):\; w(r)q'(r)\in C[0,1], \;
  w^2(r)q''(r)\in C[0,1]\},
\eqno{(4.9)}
$$
  and $\|q\|'_{C_V^1[0,1]}=\|q\|_\infty+\|wq'\|_\infty$ is an equivalent norm
  in $C_V^1[0,1]$. For $q\in {\rm Dom}(\tilde{\mathcal L}_0)$ we have, by
  definition, $\tilde{\mathcal L}_0q={\mathcal L}_0q$. Now let $f\in C^1_V[0,1]$
  and ${\rm Re}\lambda>0$. Let $q$ be the solution of (4.4). Using (4.9)
  we can easily verify that $q\in {\rm Dom}(\tilde{\mathcal L}_0)$,
  so that it is the solution of the equation $\tilde{\mathcal L}_0q-\lambda
  q=f$. Moreover, a simple computation shows that $wq'=R({\mathcal L}_0,
  \lambda)(wf')$. Thus, by (4.8) we have
$$
\begin{array}{rl}
  \|q\|'_{C^1_V[0,1]} &=\|q\|_\infty+\|wq'\|_\infty
  =\|R({\mathcal L},\lambda)f\|_\infty+
  \|R({\mathcal L},\lambda)(wf')\|_\infty\\
  &\leq \displaystyle {1\over {\rm Re}\lambda}\|f\|_\infty
  +{1\over {\rm Re}\lambda}\|wf'\|_\infty
  =\displaystyle {1\over {\rm Re}\lambda}\|f\|'_{C^1_V[0,1]}.
\end{array}
$$
  Hence $\lambda\in\rho(\tilde{\mathcal L}_0)$ and $\|R(\lambda,\tilde{\mathcal
  L}_0)\|'_{C^1_V[0,1]}\leq ({\rm Re}\lambda)^{-1}$. The desired assertion then
  follows from the Hille-Yosida Theorem and the footnote on Page 12. $\quad\Box$
\medskip

  Given $V=(\varphi,\zeta)\in S_\epsln$, we set $p(r)=p_*(r)+\varphi(r)$, $z=
  z_*+\zeta$, and as before we denote
$$
  u_{p,z}(r)={1\over r^2}\int_0^r [-K_D(c(\rho,z))+K_M(c(\rho,z))
  p(\rho)]\rho^2d\rho \quad \mbox{and} \quad
  w_{p,z}(r)=u_{p,z}(r)-ru_{p,z}(1).
$$
  Since $\|\varphi\|_\infty\leq\epsln$ and $|\zeta|\leq\epsln$, by a simple
  computation we see that
$$
  -C\epsln r(1-r)\leq w_{p,z}(r)-u_*(r)\leq C\epsln r(1-r) \quad
  \mbox{for}\;\; 0\leq r\leq 1.
\eqno{(4.10)}
$$
  Since $-C_1r(1-r)\leq u_*(r)\leq -C_2r(1-r)$, it follows that
$$
  (1+C\epsln)u_*(r)\leq w_{p,z}(r)\leq (1-C\epsln)u_*(r) \quad
  \mbox{for}\;\; 0\leq r\leq 1,
$$
  and, consequently, for $\epsln$ sufficiently small we have
$$
  -C_1r(1-r)\leq w_{p,z}(r)\leq -C_2r(1-r) \quad
  \mbox{for}\;\; 0\leq r\leq 1.
\eqno{(4.11)}
$$
  Later on we shall also use the notation $w_V(r)$ to re-denote $w_{p,z}(r)$. We
  note that all constants $C$, $C_1$ and $C_2$ that appear in (4.10)--(4.11)
  are independent of $V$ and $\epsln$.
\medskip

  {\bf Lemma 4.2}\ \ {\em $\{{\mathbb A}(V):\; V\in S_\epsln\}$ is a stable
  family of infinitesimal generators of $C_0$ semigroups on $X=C[0,1]\times
  {\mathbb R}$, and $\{\tilde{\mathbb A}(V):\; V\in S_\epsln\}$ is a stable
  family of infinitesimal generators of $C_0$ semigroups on $X_0$.}
\medskip

  {\em Proof}:\ \ Let ${\mathcal L}_Vq(r)=-w_V(r)q'(r)$. Then by Lemma 4.1 we
  know that for any $V\in S_\epsln$, ${\mathcal L}_V$ is an infinitesimal
  generator of a $C_0$ semigroup of contractions $e^{t{\mathcal L}_V}$ on
  $C[0,1]$. Since
$$
   {\mathbb A}_0(U_*+V)=\left(
\begin{array}{cc}
   {\mathcal L}_V\;\; & \;\; 0 \\ 0 \;\; & \;\; 0
\end{array}
\right),
\eqno{(4.12)}
$$
  (see (2.30)), it is evident that for any $V\in S_\epsln$, ${\mathbb A}_0
  (U_*+V)$ is an infinitesimal generator of a $C_0$ semigroup of contractions
  $e^{t{\mathbb A}_0(U_*+V)}$ on $X=C[0,1]\times {\mathbb R}$. In fact,
$$
  e^{t{\mathbb A}_0(U_*+V)}=\left(
\begin{array}{cc}
   e^{t{\mathcal L}_V}\;\; & \;\; 0\\ 0 \;\; & \;\; id
\end{array}
\right).
$$
  Hence, $\{{\mathbb A}_0(U_*+V):\; V\in S_\epsln\}$ is a stable family of
  infinitesimal generators of $C_0$ semigroups on $X$, with stability constants
  $(M,\omega)=(1,0)$. Since ${\mathbb A}(V)={\mathbb A}_0(U_*+V)+{\mathbb B}$
  and ${\mathbb B}$ is a bounded linear operator on $X$ independent of $V$, by
  a standard perturbation result (see, e. g. Theorem 2.3 in Section 5.2 of
  \cite{Pazy}) we immediately get the assertion that
  $\{{\mathbb A}(V): V\in S_\epsln\}$ is a stable family of infinitesimal
  generators of $C_0$ semigroups on $X=C[0,1]\times {\mathbb R}$, with
  stability constants $(M,\omega)=(1,\|{\mathbb B}\|)$.

  In order to prove that $\{\tilde{\mathbb A}(V):\; V\in S_\epsln\}$ is a stable
  family of infinitesimal generators of $C_0$ semigroup on $X_0=C^1_V[0,1]
  \times {\mathbb R}$, we first establish an estimate for the semigroup
  $e^{t\tilde{\mathcal L}_V}$ on $C^1_V[0,1]$ different from (4.5), where
  $\tilde{\mathcal L}_V$ represents the part of ${\mathcal L}_V$ in
  $C^1_V[0,1]$. Let $q_0\in C^1_V[0,1]$ and $q=e^{t\tilde{\mathcal L}_V}q_0=
  e^{t{\mathcal L}_V}q_0$. Then $q$ is the solution of the problem:
$$
  {\partial q\over\partial t}+w_V(r){\partial q\over\partial r}=0 \quad
  \mbox{for}\;\; 0\leq r\leq 1\;\;\mbox{and} \;\; t>0, \quad
  q|_{t=0}=q_0.
$$
  Let $l(r,t)=r(1\!-\!r)\displaystyle{\partial q(r,t)\over\partial r}$ and
  $l_0(r)=r(1\!-\!r)q_0'(r)$. Differentiating the above equation in $r$ and
  multiplying it with $r(1\!-\!r)$, we get
$$
  {\partial l\over\partial t}+w_V(r){\partial l\over\partial r}=a_V(r)l
  \quad \mbox{for}\;\; 0\leq r\leq 1\;\;\mbox{and} \;\; t>0, \quad
  l|_{t=0}=l_0,
$$
  where $a_V(r)=(1-2r)\displaystyle{w_V(r)\over r(1\!-\!r)}-w_V'(r)$. Clearly,
  there exists a nonnegative constant $c_0$ independent of $V$ such that
$$
  a_V(r)\leq c_0 \quad \mbox{for}\;\; 0<r<1,\quad
  \mbox{for all}\;\; V\in S_\epsln.
$$
  Using this fact and a standard characteristics argument we can easily obtain
$$
  \|l(\cdot,t)\|_\infty\leq\|l_0\|_\infty e^{c_0 t} \quad
  \mbox{for}\;\; t\geq 0.
$$
  Combining this estimate with $\|q(\cdot,t)\|_\infty\leq\|q_0\|_\infty$
  ensured by (4.2) we get
$$
  \|q(\cdot,t)\|_{C^1_V[0,1]}\leq\|q_0\|_{C^1_V[0,1]} e^{c_0 t} \quad
  \mbox{for}\;\; t\geq 0.
$$
  Hence
$$
  \|e^{t\tilde{\mathcal L}_V}\|_{L(C^1_V[0,1])}\leq e^{c_0 t} \quad
  \mbox{for}\;\; t\geq 0,\quad \mbox{for all}\;\; V\in S_\epsln.
$$
  Hence, $\{\tilde{\mathcal L}_V:V\in S_\epsln\}$ is a stable family of
  infinitesimal generators of $C_0$ semigroups on $C^1_V[0,1]$, with stability
  constants $(M,\omega)=(1,c_0)$. Using this assertion and (4.12) we see
  easily that $\{\tilde{\mathbb A}_0(U_*+V): V\in S_\epsln\}$, the part of
  $\{{\mathbb A}_0(U_*+V): V\in S_\epsln\}$ on $X_0=C[0,1]\times {\mathbb R}$,
  is a stable family of infinitesimal generators of $C_0$ semigroups on $X_0$,
  with stability constants $(M,\omega)=(1,c_0)$. Since $\tilde{\mathbb A}(V)
  =\tilde{\mathbb A}_0(U_*+V)+{\mathbb B}$ and, by Corollary 3.2, ${\mathbb B}$
  is a bounded linear operator on $X_0$ independent of $V$, we conclude as
  before that $\{{\mathbb A}(V): V\in S_\epsln\}$ is a stable family of
  infinitesimal generators of $C_0$ semigroups on $X_0=C^1_V[0,1]\times
  {\mathbb R}$, with stability constants $(M,\omega)=
  (1,c_0+\|{\mathbb B}\|_{L(C^1_V[0,1])})$. This completes the proof of Lemma
  4.2. $\quad\Box$
\medskip

  Since ${\mathbb A}\in C^\infty(X,L(X_0,X))$, by Lemma 4.2 we see that for any
  $V\in C([0,\infty),X)$ such
  that $V(t)\in S_\epsln$ for all $t\geq 0$, $\{{\mathbb A}(V(t)):t\geq 0\}$
  satisfies the conditions $(H_1)$--$(H_3)$ in Section 5.3 of \cite{Pazy}. It
  follows by Theorem 3.1 in Section 5.3 of \cite{Pazy} that given a such
  function $V=V(t)$, there exists an evolution system determined by
  $\{{\mathbb A}(V(t)):t\geq 0\}$, which we denote as ${\mathbb U}(t,s,V)$. By
  definition, this means that

  (1) for any $t\geq s\geq 0$, ${\mathbb U}(t,s,V)$ is a bounded linear
  operator on $X$,

  (2) ${\mathbb U}(s,s,V)=id$ for all $s\geq 0$, ${\mathbb U}(t,s,V){\mathbb U}
  (s,r,V)={\mathbb U}(t,r,V)$ for all $t\geq s\geq r$, and

  (3) the mapping $(t,s)\to {\mathbb U}(t,s,V)$ is strongly continuous for
  $t\geq s\geq 0$.

\noindent
  However, the theory developed in \cite{Pazy} does not ensure that $U=
  {\mathbb U}(t,s,V)U_0$ is a solution of the problem
$$
\left\{
\begin{array}{l}
  \displaystyle{dU\over dt}={\mathbb A}(V(t))U \;\;\; \mbox{for}\;\; t>s,
  \\ [0.1cm]
  U|_{t=s}=U_0,
\end{array}
\right.
\eqno{(4.13)}
$$
  even if $U_0\in X_0$, unless some other conditions are satisfied by
  ${\mathbb U}(t,s,V)$. These conditions are as follows (see the conditions
  $(E_4)$ and $(E_5)$ in Theorem 4.3 in Section 5.4 of \cite{Pazy}):

  (4) ${\mathbb U}(t,s,V)X_0\subseteq X_0$ for any $t\geq s\geq 0$, and

  (5) for any $U_0\in X_0$, the mapping $(t,s)\to {\mathbb U}(t,s,V)U_0$ is
  continuous in $X_0$ for $t\geq s\geq 0$.

\noindent
  In the following lemma we shall directly prove that for any $U_0\in X_0$,
  the problem (4.13) has a unique solution $U=U_s(t)\in C([s,\infty),X_0)\cap
  C^1([s,\infty),X)$. By Theorem 4.2 in Section 5.4 of \cite{Pazy}, it then
  follows that $U_s(t)={\mathbb U}(t,s,V)U_0$ and, consequently, the conditions
  (4) and (5) above are satisfied.
\medskip

  {\bf Lemma 4.3}\ \ {\em Given $V\in C([0,\infty),X)$ such that $V(t)\in
  S_\epsln$ for all $t\geq 0$, for any $s\geq 0$ and any $U_0\in X_0$ the
  problem $(4.13)$ has a unique solution $U=U_s(t)\in C([s,\infty),X_0)\cap
  C^1([s,\infty),X)$.}
\medskip

  {\em Proof}:\ \ Let $U=(q,y)$ and $U_0=(q_0,y_0)$. Then (4.13) can be
  rewritten as follows:
$$
\left\{
\begin{array}{l}
  \displaystyle {\partial q\over\partial t}+w_V(r,t){\partial q\over\partial r}
  =a(r)q+{\mathcal B}(q)+b(r)y
  \quad \mbox{for} \;\; 0\leq r\leq 1, \;\; t>s,\\ [0.1cm]
  \displaystyle {dy\over dt}={\mathcal F}(q)+\kappa y
  \quad \mbox{for} \;\; t>s,\\ [0.1cm]
  q|_{t=s}=q_0(r)\quad \mbox{for} \;\; 0\leq r\leq 1, \quad
  \mbox{and}\quad y|_{t=s}=y_0.
\end{array}
\right.
\eqno{(4.14)}
$$
  Using the characteristic method and the Banach fixed point theorem, we can
  easily show that this problem has a unique local solution $(q,y)$ with $q\in
  C([0,1]\times [0,\delta])$ and $y\in C^1[0,\delta]$ for some $\delta>0$.
  Since $w_V(0,t)=w_V(1,t)=0$ for all $t\geq 0$,  we see that the two lines
  $r=0$ and $r=1$ are characteristic curves. It follows that all characteristic
  curves starting from the open interval $(0,1)$ always lie in it, so that the
  solution of the above problem exists for all $t\geq s$. It remains to prove
  that $q\in C([0,\infty),C^1_V[0,1])$. To this end we formally differentiate
  the first equation in (4.14) in $r$ and multiply it with $r(1-r)$, which gives,
  by letting $l(r,t)=r(1\!-\!r)\displaystyle{\partial q(r,t)\over\partial r}$,
  that
$$
  {\partial l\over\partial t}+w_V(r,t){\partial l\over\partial r}
  =a_1(r,t)u+f_1(r,t)\quad \mbox{for} \;\; 0\leq r\leq 1, \;\; t>0,
\eqno{(4.15)}
$$
  where
$$
  a_1(r,t)=a(r)+(1-2r){w_V(r,t)\over r(1\!-\!r)}-
  {\partial w_V(r,t)\over\partial r},
$$
$$
  f_1(r,t)=r(1\!-\!r)a'(r)q(r,t)+r(1\!-\!r)
  {\partial{\mathcal B}q(r,t)\over\partial r}+r(1\!-\!r)b'(r)y(t).
$$
  Clearly, $a_1\in C([0,1]\times [0,\infty))$. By Corollary 3.2 we see that also
  $f_1\in C([0,1]\times [0,\infty)$. Thus by using the characteristic method we
  can easily prove that (4.15) imposed with the initial condition $l(r,0)=
  r(1\!-\!r)q_0'(r)$ has a unique solution $l\in C([0,1]\times [0,\infty)$.
  Thus, the above formal computation makes sense and, consequently, $q\in
  C([0,\infty),C^1_V[0,1])\cap C^1([0,\infty),C[0,1])$. The desired assertion
  now becomes immediate. $\quad\Box$
\medskip

  By the above results and Theorems 4.2 and 5.2 in Sections 5.4 and 5.5 of
  \cite{Pazy}, we get:
\medskip

  {\bf Corollary 4.4}\ \ {\em Let $V=V(t)\in C([0,\infty),X)$ be as in Lemma
  4.3, and let $F=F(t)\in C([0,\infty),X_0)$. Then for any $U_0\in X_0$, the
  initial value problem
$$
  {dU\over dt}={\mathbb A}(V(t))U+F(t) \quad \mbox{for}\;\; t>0, \quad
  U(0)=U_0
$$
  has a unique solution $U=U(t)\in C([0,\infty),X_0)\cap C^1([0,\infty),X)$,
  and it is given by
$$
  U(t)={\mathbb U}(t,0,V)U_0+\int_0^t {\mathbb U}(t,s,V)F(s) ds.
$$
  $\Box$}

\section{Similarity transformation}

  In this section we shall study a family of $C^1$-diffeomorphisms $\bar{r}=
  T(r,t,s)$ of the unit interval $0\leq r\leq 1$ to itself, where $t\geq s\geq
  0$ are parameters. This family of diffeomorphisms will be used in the next
  section to deduce a uniform decay estimate for the evolution system ${\mathbb
  U}(t,s, V)$ established in the previous section when $V$ is replaced by an
  exponentially decaying function $V=V(t)\in C([0,\infty),S_\epsln)$.

  Let $w\in C([0,\infty),C^1[0,1])$. We assume that $w$ satisfies the following
  condition: For some small parameter $\epsln>0$,
$$
  -C\epsln r(1-r)e^{-\mu t}\leq w(r,t)-u_*(r)\leq C\epsln r(1-r)e^{-\mu t}
  \quad \mbox{for}\;\; 0\leq r\leq 1, \;\; t\geq 0,
\eqno{(5.1)}
$$
  where $C$ is a positive constant independent of $\epsln$ and $w$.
  Since $-C_1r(1-r)\leq u_*(r)\leq -C_2r(1-r)$, we see that
$$
  \sup_{0<r<1}\Big|{w(r,t)\over u_*(r)}-1\Big|\leq C\epsln e^{-\mu t}
  \quad \mbox{for}\;\; 0\leq r\leq 1, \;\; t\geq 0,
\eqno{(5.2)}
$$
  and for $\epsln$ sufficiently small we have
$$
  -C_1r(1-r)\leq w(r,t)\leq -C_2r(1-r) \quad
  \mbox{for}\;\; 0\leq r\leq 1, \;\; t\geq 0
\eqno{(5.3)}
$$
  and
$$
  {1\over 2}\leq {w(r,t)\over u_*(r)}\leq 2
  \quad \mbox{for}\;\; 0\leq r\leq 1, \;\; t\geq 0.
\eqno{(5.4)}
$$
\medskip

  Let $0\leq\xi\leq 1$ and $s\geq 0$. Consider the following initial value
  problem:
$$
  \displaystyle {dr\over dt}=u_*(r) \quad \mbox{for}\;\; t>s, \quad
  \displaystyle r|_{t=s}=\xi.
\eqno{(5.5)}
$$
  Since $u_*\in C^1[0,1]$, $u_*(r)<0$ for $0<r<1$ and, in particular, $u_*(0)=
  u_*(1)=0$, it can be easily shown that this problem has a unique solution
  $r=\Phi_*(\xi,t,s)$ for all $t\geq s$, satisfying the following properties:
$$
  \Phi_*(\xi,t,s)\;\; \mbox{is twice continuously differentiable in}\;\;
  (\xi,t,s),
$$
$$
  \Phi_*(0,t,s)=0, \quad \Phi_*(1,t,s)=1 \quad \mbox{for}\;\; t\geq s,
$$
$$
  0<\Phi_*(\xi,t,s)<1 \quad \mbox{for}\;\; 0<\xi<1,\;\; t\geq s,
$$
$$
  {\partial\Phi_*(\xi,t,s)\over\partial\xi}>0, \quad
  {\partial\Phi_*(\xi,t,s)\over\partial t}<0 \quad
  \mbox{for}\;\; 0<\xi<1,\;\; t\geq s.
$$
  Note that we also have
$$
  \Phi_*(\xi,s,s)=\xi \quad \mbox{and} \quad
  \Phi_*(\xi,t,s)=\Phi_*(\xi,t-s,0) \quad
  \mbox{for}\;\; 0\leq\xi\leq 1,\;\; t\geq s.
$$
  From these properties we see that for any $s\geq 0$ and $t\geq s$,
  the mapping $\xi\to r=\Phi_*(\xi,t,s)$ is a $C^2$ diffeomorphism of $[0,1]$
  to itself. Let $\xi=\Psi_*(r,t,s)$ be the inverse of this mapping. Clearly,
  $\Psi_*$ satisfies the following properties:
$$
  \Psi_*(r,t,s)\;\; \mbox{is twice continuously differentiable in}\;\;
  (r,t,s),
$$
$$
  \Psi_*(0,t,s)=0, \quad \Psi_*(1,t,s)=1 \quad \mbox{for}\;\; t\geq s,
$$
$$
  0<\Psi_*(r,t,s)<1 \quad \mbox{for}\;\; 0<r<1,\;\; t\geq s,
$$
$$
  {\partial\Psi_*(r,t,s)\over\partial r}>0, \quad
  {\partial\Psi_*(r,t,s)\over\partial t}>0 \quad
  \mbox{for}\;\; 0<r<1,\;\; t\geq s,
$$
$$
  \Psi_*(r,s,s)=r \quad \mbox{and} \quad \Psi_*(r,t,s)=\Psi_*(r,t-s,0)
  \quad \mbox{for}\;\; 0\leq r\leq 1,\;\; s\geq 0.
$$
  Furthermore, by the definition of $\Psi_*$ we have the following relations:
$$
  \Psi_*(\Phi_*(\xi,t,s),t,s)=\xi \quad \mbox{for} \;\;
  0\leq\xi\leq 1,\;\; t\geq s,
$$
$$
  \Phi_*(\Psi_*(r,t,s),t,s)=r \quad \mbox{for} \;\;
  0\leq r\leq 1,\;\; t\geq s.
$$
  From the first relation we easily deduce that $\xi=\Psi_*(r,t,s)$ is the
  unique solution of the following initial value problem:
$$
  \displaystyle {\partial\xi\over\partial t}+
  u_*(r){\partial\xi\over\partial r}=0 \quad \mbox{for}\;\; t>s,\quad
  \displaystyle \xi|_{t=s}=r.
\eqno{(5.6)}
$$

  Next, let $r=\Phi(\xi,t,s)$ ($0\leq\xi\leq 1$, $t\geq s\geq 0$) be the
  solution of the following problem:
$$
  \displaystyle {dr\over dt}=w(r,t) \quad \mbox{for}\;\; t>s,\quad
  \displaystyle r|_{t=s}=\xi.
\eqno{(5.7)}
$$
  Similarly as before, $\Phi(\xi,t,s)$ is well-defined for all $0\leq\xi\leq
  1$ and $t\geq s$, and it satisfies the following properties:
$$
  \Phi(\xi,t,s)\;\; \mbox{is continuously differentiable in}\;\; (\xi,t,s),
$$
$$
  \Phi(0,t,s)=0, \quad \Phi(1,t,s)=1 \quad \mbox{for}\;\; t\geq s,
$$
$$
  0<\Phi(\xi,t,s)<1 \quad \mbox{for}\;\; 0<\xi<1,\;\; t\geq s,
$$
$$
  {\partial\Phi(\xi,t,s)\over\partial\xi}>0, \quad
  {\partial\Phi(\xi,t,s)\over\partial t}<0 \quad
  \mbox{for}\;\; 0<\xi<1,\;\; t\geq s,
$$
$$
  \Phi(\xi,s,s)=\xi \quad \mbox{for}\;\; 0\leq\xi\leq 1,\;\; s\geq 0.
$$
  From the above properties we see that for any $s\geq 0$ and $t\geq s$, the
  mapping $\xi\to r=\Phi(\xi,t,s)$ is a $C^1$ diffeomorphism of $[0,1]$ to
  itself. Let $\xi=\Psi(r,t,s)$ be the inverse of this mapping. Similarly as
  before we have
$$
  \Psi(r,t,s)\;\; \mbox{is continuously differentiable in}\;\;
  (r,t,s),
$$
$$
  \Psi(0,t,s)=0, \quad \Psi(1,t,s)=1 \quad \mbox{for}\;\; t\geq s,
$$
$$
  0<\Psi(r,t,s)<1 \quad \mbox{for}\;\; 0<r<1,\;\; t\geq s,
$$
$$
  {\partial\Psi(r,t,s)\over\partial r}>0, \quad
  {\partial\Psi(r,t,s)\over\partial t}>0 \quad
  \mbox{for}\;\; 0<r<1,\;\; t\geq s,
$$
$$
  \Psi(r,s,s)=r \quad \mbox{for}\;\; 0\leq r\leq 1,\;\; s\geq 0.
$$
  Moreover, we have the following relations:
$$
  \Psi(\Phi(\xi,t,s),t,s)=\xi \quad \mbox{for} \;\;
  0\leq\xi\leq 1,\;\; t\geq s,
$$
$$
  \Phi(\Psi(r,t,s),t,s)=r \quad \mbox{for} \;\;
  0\leq r\leq 1,\;\; t\geq s,
$$
  and $\xi=\Psi(r,t,s)$ is the unique solution of the following initial value
  problem:
$$
  \displaystyle {\partial\xi\over\partial t}+
  w(r,t){\partial\xi\over\partial r}=0 \quad \mbox{for}\;\; t>s,\quad
  \displaystyle \xi|_{t=s}=r.
\eqno{(5.8)}
$$

  In the sequel we consider the following initial value problem:
$$
\left\{
\begin{array}{l}
  \displaystyle {\partial \bar{r}\over\partial t}+w(r,t){\partial \bar{r}\over\partial r}
  =u_*(\bar{r}) \quad \mbox{for} \;\; 0\leq r\leq 1,\;\; t>s,\\[0.2cm]
  \displaystyle \bar{r}|_{t=s}=r \quad \mbox{for} \;\; 0\leq r\leq 1.
\end{array}
\right.
\eqno{(5.9)}
$$

  {\bf Lemma 5.1}\ \ {\em For any $0\leq r\leq 1$ and $s\geq 0$, the problem
  $(5.9)$ has a unique solution $\bar{r}=T(r,t,s)$ for all $t\geq s$, and the
  following relation holds:
$$
  T(r,t,s)=\Phi_*(\Psi(r,t,s),t,s) \quad \mbox{for}\;\;
  0\leq r\leq 1,\;\; t\geq s\geq 0.
\eqno{(5.10)}
$$
}

  {\em Proof}:\ \ Using (5.5) and (5.8) we can easily verify that $\bar{r}=
  \Phi_*(\Psi(r,t,s),t,s)$ is a solution of the problem (5.9). Thus, (5.10)
  follows by uniqueness of the solution. $\quad\Box$
\medskip

  By (5.10), it is evident that for any $s\geq 0$ and $t\geq s$, the mapping
  $r\to \bar{r}=T(r,t,s)$ is a $C^1$ diffeomorphism of $[0,1]$ to itself, satisfying
  the following properties:
$$
  T(0,t,s)=0, \quad T(1,t,s)=1 \quad \mbox{for}\;\; t\geq s\geq 0,
$$
$$
  {\partial T(r,t,s)\over\partial r}>0 \quad
  \mbox{for}\;\; 0<r<1,\;\; t\geq s.
$$
  We denote by $r=S(\bar{r},t,s)$ the inverse of this mapping. By (5.10) it is clear
  that
$$
  S(\bar{r},t,s)=\Phi(\Psi_*(\bar{r},t,s),t,s) \quad
  \mbox{for}\;\; 0\leq \bar{r}\leq 1,\;\; t\geq s\geq 0.
\eqno{(5.11)}
$$
  It is also clear that $S(\bar{r},t,s)$ satisfies the following properties:
$$
  S(0,t,s)=0, \quad S(1,t,s)=1 \quad \mbox{for}\;\; t\geq s\geq 0,
$$
$$
  {\partial S(\bar{r},t,s)\over\partial \bar{r}}>0 \quad
  \mbox{for}\;\; 0<\bar{r}<1,\;\; t\geq s\geq 0.
$$

  $T$ and $S$ can be expressed in more explicit formulations. To show this we
  introduce a function $F_*$ as follows:
$$
  F_*(r)=-\int_{{1\over 2}}^r{d\eta\over u_*(\eta)}
  =\int_{{1\over 2}}^r{d\eta\over |u_*(\eta)|} \quad
  \mbox{for}\;\; 0<r<1.
\eqno{(5.12)}
$$
  Clearly, $F_*\in C^1(0,1)$, $F_*'(r)>0$ for all $0<r<1$, and
$$
  \lim_{r\to 0^+}F_*(r)=-\infty, \quad
  \lim_{r\to 1^-}F_*(r)=\infty.
$$
  Hence $\bar{r}=F_*(r)$ is a $C^1$ diffeomorphism of the open unit interval
  $(0,1)$ to the real line $(-\infty,\infty)$. From (5.5) we easily obtain
$$
  F_*(\Phi_*(\xi,t,s))-F_*(\xi)=-t+s.
$$
  Thus
$$
  \Phi_*(\xi,t,s)=F_*^{-1}(F_*(\xi)-t+s),
\eqno{(5.13)}
$$
  and, consequently,
$$
  \Psi_*(r,t,s)=F_*^{-1}(F_*(r)+t-s).
\eqno{(5.14)}
$$
  Next, let
$$
  g(\xi,t,s)=G(\Phi(\xi,t,s),t),\quad \mbox{where}\;\;\;
  G(r,t)={w(r,t)\over u_*(r)}-1.
$$
  Since $w(r,t)=[1+G(r,t)]u_*(r)$, from (5.7) we see that $r=\Phi(\xi,t,s)$
  is a solution of the following problem:
$$
  \displaystyle {dr\over dt}=[1+g(\xi,t,s)]u_*(r) \quad
  \mbox{for}\;\; t>s,\quad \displaystyle r|_{t=s}=\xi.
\eqno{(5.15)}
$$
  Thus similarly as before we have
$$
  F_*(\Phi(\xi,t,s))-F_*(\xi)=-t+s-\int_s^t g(\xi,\tau,s)d\tau,
\eqno{(5.16)}
$$
  so that
$$
  \Phi(\xi,t,s)=F_*^{-1}\big(F_*(\xi)-t+s-\int_s^t g(\xi,\tau,s)d\tau\big),
\eqno{(5.17)}
$$
$$
  \Psi(r,t,s)=F_*^{-1}\big(F_*(r)+t-s+\int_s^t g(\Psi(r,t,s),\tau,s)d\tau\big).
\eqno{(5.18)}
$$
  Combining (5.10), (5.11), (5.13), (5.14), (5.17) and (5.18) we see that
$$
  T(r,t,s)=F_*^{-1}(F_*(r)+\int_s^t g(\Psi(r,t,s),\tau,s)d\tau),
\eqno{(5.19)}
$$
$$
  S(\bar{r},t,s)=F_*^{-1}(F_*(\bar{r})-\int_s^t
  g(\Psi_*(\bar{r},t,s),\tau,s)d\tau).
\eqno{(5.20)}
$$
\medskip

  {\bf Lemma 5.2}\ \ {\em  Assume that $|\zeta|\leq C$. Then there exist
  positive constants $C_1$ and $C_2$ depending only on $C$ such that for any
  $0<r<1$ we have
$$
  C_1 r(1-r)\leq F_*^{-1}(F_*(r)+\zeta)\big[1-F_*^{-1}(F_*(r)+\zeta)\big]
  \leq C_2 r(1-r).
\eqno{(5.21)}
$$
}

  {\em Proof}:\ \ Since $-C\leq\zeta\leq C$, by the monotonicity of $F_*$
  we have
$$
  F_*^{-1}(F_*(r)-C)\leq F_*^{-1}(F_*(r)+\zeta)\leq F_*^{-1}(F_*(r)+C).
$$
  Thus
$$
  {F_*^{-1}(F_*(r)+\zeta)\over r}\leq {F_*^{-1}(F_*(r)+C)\over r}
$$
  and
$$
  {1-F_*^{-1}(F_*(r)+\zeta)\over 1-r}\leq {1-F_*^{-1}(F_*(r)-C)\over 1-r}.
$$
  We claim that
$$
  \lim_{r\to 0^+}{F_*^{-1}(F_*(r)+C)\over r}=e^{C|u_*'(0)|}
\eqno{(5.22)}
$$
  Indeed, since $u_*(r)=u_*'(0)r[1+O(r^\beta)]$ (for $r\sim 0$) for some $0<
  \beta\leq 1$ (see Assertion (3) of Lemma 3.1), we have $1/u_*(r)=
  [1+O(r^\beta)]/u_*'(0)r$ (for $r\sim 0$), so that
$$
  F_*(r)=-\int_{{1\over 2}}^r {d\eta\over u_*(\eta)}=
  -\int_{r_0}^r {d\eta\over u_*(\eta)}
  -\int_{{1\over 2}}^{r_0}{d\eta\over u_*(\eta)}=
  -\int_{r_0}^r {1+O(\eta^\beta)\over u_*'(0)\eta}d\eta+C,
$$
  which yields
$$
  F_*(r)=\log r^{{1\over |u_*'(0)|}}+C_1+O(r^\beta) \quad
  \mbox{for}\;\; r\sim 0.
$$
  Thus
$$
  F_*^{-1}(\xi)=e^{|u_*'(0)|[\xi-C_1+O\big((F_*^{-1}(\xi))^\beta\big)]}
  \quad \mbox{for}\;\; \xi\sim -\infty,
$$
  and, consequently,
$$
  F_*^{-1}(F_*(r)+C)=re^{C|u_*'(0)|+O(r^\beta)} \quad \mbox{for}\;\; r\sim 0,
$$
  by which (5.22) follows immediately. Similarly, we also have
$$
  \lim_{r\to 1^-}{1-F_*^{-1}(F_*(r)-C)\over 1-r}=e^{C u_*'(1)}.
\eqno{(5.23)}
$$
  By (5.22) and (5.23), the second inequality in (5.21) immediately follows.
  The proof for the first inequality in (5.21) is similar. This
  completes the proof of Lemma 5.2. $\quad\Box$
\medskip

  {\bf Corollary 5.3}\ \ {\em For $\epsln$ sufficiently small we have
$$
  C_1 r(1-r)\leq T(r,t,s)[1-T(r,t,s)]\leq C_2 r(1-r),
\eqno{(5.24)}
$$
$$
  C_1 \bar{r}(1-\bar{r})\leq S(\bar{r},t,s)[1-S(\bar{r},t,s)]\leq C_2
  \bar{r}(1-\bar{r}).
\eqno{(5.25)}
$$
}

  {\em Proof}:\ \ Let $\zeta=\displaystyle\int_s^t g(\Psi(r,t,s),\tau,s)d\tau$.
  By (5.2) we have
$$
  |\zeta|\leq\int_s^t\!\!|g(\Psi(r,t,s),\tau,s)|d\tau\leq
  C\epsln\int_s^t e^{-\mu\tau}d\tau\leq C\epsln\int_0^\infty
  e^{-\mu\tau}d\tau\leq C\epsln\leq C.
$$
  Hence, (5.24) follows from (5.19) and (5.21). Similarly, (5.25) follows from
  (5.20) and (5.21). $\quad\Box$
\medskip

  As an immediate consequence of Corollary 5.3 we see that there exists constant
  $C>1$ such that for $\epsln$ sufficiently small we have
$$
  C^{-1}\leq {T(r,t,s)\over r}\leq C \quad \mbox{and} \quad
  C^{-1}\leq {S(\bar{r},t,s)\over \bar{r}}\leq C.
$$

  {\bf Corollary 5.4}\ \ {\em  For $\epsln$ sufficiently small we have the
  following inequalities:
$$
  C_1\Psi_*(r,t,s)[1-\Psi_*(r,t,s)]\leq\Psi(r,t,s)[1-\Psi(r,t,s)]\leq
  C_2\Psi_*(r,t,s)[1-\Psi_*(r,t,s)].
\eqno{(5.26)}
$$
$$
  C_1\Phi_*(\bar{r},t,s)[1-\Phi_*(\bar{r},t,s)]\leq\Phi(\bar{r},t,s)
  [1-\Phi(\bar{r},t,s)]\leq C_2\Phi_*(\bar{r},t,s)[1-\Phi_*(\bar{r},t,s)].
\eqno{(5.27)}
$$
}

  {\em Proof}:\ \ Let $\bar{r}=\Psi_*(r,t,s)$ and $\zeta\!=\!\displaystyle\int_s^t
  \!\!g(\Psi(r,t,s),\tau,s)d\tau$. Then by (5.14) and (5.18) we have
$$
  \Psi(r,t,s)=F_*^{-1}(F_*(\bar{r})+\zeta).
$$
  By this expression and (5.21) we immediately obtain (5.26). The proof of
  (5.27) is similar. $\quad\Box$
\medskip

  As an immediate consequence of Corollary 5.4 we see that there exists constant
  $C>1$ such that for $\epsln$ sufficiently small we have
$$
  C^{-1}\leq {\Psi(r,t,s)\over\Psi_*(r,t,s)}\leq C \quad \mbox{and} \quad
  C^{-1}\leq {\Phi(\bar{r},t,s)\over\Phi_*(\bar{r},t,s)}\leq C.
$$

  {\bf Lemma 5.5}\ \ {\em  We have the following inequalities:
$$
  |T(r,t,s)-r|\leq C\epsln\big(e^{-\mu s}-e^{-\mu t}\big)r(1-r),
\eqno{(5.28)}
$$
$$
  |S(\bar{r},t,s)-\bar{r}|\leq C\epsln\big(e^{-\mu s}-e^{-\mu t}\big)
  \bar{r}(1-\bar{r}).
\eqno{(5.29)}
$$
}

  {\em Proof}:\ \ Similarly as in the proof of Corollary 5.3 we have
$$
  \int_s^t |g(\Psi(r,t,s),\tau,s)|d\tau\leq C\epsln\int_s^t e^{-\mu\tau}
  \leq C\epsln\Big(e^{-\mu s}-e^{-\mu t}\Big).
$$
  Thus, by noticing that $\displaystyle{dF_*^{-1}(\eta)\over d\eta}=
  {1\over F_*'(F_*^{-1}(\eta))}=|u_*(F_*^{-1}(\eta))|$, we see that
$$
\begin{array}{rl}
  |T(r,t,s)-r|&=\displaystyle|F_*^{-1}(F_*(r)+\int_s^t g(\Psi(r,t,s),
  \tau,s)d\tau)-F_*^{-1}(F_*(r))|\\ [0.3cm]
  &\leq\displaystyle\int_0^1|u_*(F_*^{-1}(F_*(r)+\zeta_\theta)))|d\theta\cdot
  \int_s^t |g(\Psi(r,t,s),\tau,s)|d\tau\\ [0.3cm]
  &\leq\displaystyle C\epsln\Big(e^{-\mu s}-e^{-\mu t}\Big)\cdot
  \int_0^1|u_*(F_*^{-1}(F_*(r)+\zeta_\theta)))|d\theta.
\end{array}
$$
  where $\zeta_\theta=\displaystyle\theta\!\int_s^t\!g(\Psi(r,t,s),\tau,s)
  d\tau$. Since $|\zeta_\theta|\leq C$ and $|u_*(\eta)|\leq C\eta(1-\eta)$, by
  Lemma 5.2 we have
$$
  |u_*(F_*^{-1}(F_*(r)+\zeta_\theta)))|\leq CF_*^{-1}(F_*(r)+\zeta_\theta))
  [1-F_*^{-1}(F_*(r)+\zeta_\theta))]\leq Cr(1-r).
$$
  Substituting this estimate into the above inequality, we see that $(5.28)$
  follows. The proof of (5.29) is similar. $\quad\Box$
\medskip

  {\bf Corollary 5.6}\ \ {\em  We have the following inequalities:
$$
  |\Phi(\xi,t,s)-\Phi_*(\xi,t,s)|\leq C\epsln\big(e^{-\mu s}-e^{-\mu t}\big)
  \Phi_*(\xi,t,s)[1-\Phi_*(\xi,t,s)],
\eqno{(5.30)}
$$
$$
  |\Psi(r,t,s)-\Psi_*(r,t,s)|\leq C\epsln\big(e^{-\mu s}-e^{-\mu t}\big)
  \Psi_*(r,t,s)[1-\Psi_*(r,t,s)].
\eqno{(5.31)}
$$
}

  {\em Proof}:\ \ Let $r=\Phi(\xi,t,s)$. Then $\xi=\Psi(r,t,s)$ and
  $\Phi_*(\xi,t,s)=\Phi_*(\Psi(r,t,s),t,s)=T(r,t,s)$. Thus by (5.28) we have
$$
  |\Phi(\xi,t,s)-\Phi_*(\xi,t,s)|=|r-T(r,t,s)|\leq C\epsln
  \big(e^{-\mu s}-e^{-\mu t}\big)r(1-r).
$$
  Substituting $r=\Phi(\xi,t,s)$ into the right-hand side of the last
  inequality and using (5.27), we see that (5.30) follows. The proof of (5.31)
  is similar. $\quad\Box$
\medskip

  {\bf Lemma 5.7}\ \ {\em Assume that in addition to $(5.2)$ there also holds
$$
  \max_{0\leq r\leq 1}\Big|{\partial w(r,t)\over\partial r}-u_*'(r)\Big|
  \leq C\epsln e^{-\mu t}.
\eqno{(5.32)}
$$
  Then we have the following estimates:
$$
  e^{-C\epsln\big(e^{-\mu s}-e^{-\mu t}\big)}\leq
  {\partial T(r,t,s)\over\partial r}\leq
  e^{C\epsln\big(e^{-\mu s}-e^{-\mu t}\big)}.
\eqno{(5.33)}
$$
$$
  e^{-C\epsln\big(e^{-\mu s}-e^{-\mu t}\big)}\leq
  {\partial S(\bar{r},t,s)\over\partial \bar{r}}\leq
  e^{C\epsln\big(e^{-\mu s}-e^{-\mu t}\big)}.
\eqno{(5.34)}
$$
}

  {\em Proof}: \ \ Recalling that $T(r,t,s)=\Phi_*(\Psi(r,t,s),t,s)$, we see
  that
$$
\begin{array}{rcl}
   \qquad\qquad \displaystyle{\partial T(r,t,s)\over\partial r}
  &=&\displaystyle{\partial\Phi_*\over\partial\xi}(\Psi(r,t,s),t,s)
    {\partial\Psi(r,t,s)\over\partial r}
   =\displaystyle{\partial\Phi_*\over\partial\xi}(\Psi(r,t,s),t,s)
    \Big[{\partial\Phi\over\partial\xi}(\Psi(r,t,s),t,s)\Big]^{-1}
\\[0.3cm]
  &=&\displaystyle \exp\Big({\int_s^t\big[u_*'(\Phi_*(\xi,\tau,s))-
    {\partial w\over\partial r}(\Phi(\xi,\tau,s),\tau)\big]d\tau}\Big)
    \Big|_{\xi=\Psi(r,t,s)}.\hfill (5.35)
\end{array}
%\eqno{(5.35)}
$$
  We have
$$
\begin{array}{rcl}
  \displaystyle u_*'(\Phi_*(\xi,\tau,s))-{\partial w\over\partial r}
  (\Phi(\xi,\tau,s),\tau)&=&\displaystyle [u_*'(\Phi_*(\xi,\tau,s))-
  u_*'(\Phi(\xi,\tau,s))]\\[0.1cm]
  && \displaystyle+[u_*'(\Phi(\xi,\tau,s))-{\partial w\over\partial r}
  (\Phi(\xi,\tau,s),\tau)].
\end{array}
$$
  By the assumption $(5.32)$ we have
$$
  \sup_{\xi\in{\mathbb R}}\Big|u_*'(\Phi(\xi,\tau,s))-{\partial w\over
  \partial r}(\Phi(\xi,\tau,s),\tau)\Big|\leq
  \sup_{0<r<1}\Big|u_*'(r)-{\partial w\over\partial r}(r,\tau)\Big|\leq
  C\epsln e^{-\mu\tau}.
$$
  Next, by the assertion (3) of Lemma 3.1 we know that there exists $0\leq
  \gamma<1$ such that $r^\gamma u_*''(r)\in C[0,1]$. With this fact in mind,
  we use the mean value theorem to compute
$$
\begin{array}{rcl}
  &&\displaystyle \Big|u_*'(\Phi_*(\xi,\tau,s))-u_*'(\Phi(\xi,\tau,s))\Big|
    \Big|_{\xi=\Psi(r,t,s)}
   =\displaystyle |u_*''(\zeta)|\big|\Phi_*(\xi,\tau,s)-\Phi(\xi,\tau,s)\big|
    \big|_{\xi=\Psi(r,t,s)} \\ [0.2cm]
  &=&\displaystyle |\zeta^\gamma u_*''(\zeta)|\cdot
    \Big[\Big({\Phi(\xi,\tau,s)\over\zeta}\Big)^\gamma\cdot
    \big(\Phi(\xi,\tau,s)\big)^{-\gamma}\big|\Phi_*(\xi,\tau,s)-
    \Phi(\xi,\tau,s)\big|\Big]\big|_{\xi=\Psi(r,t,s)},
\end{array}
$$
  where $\zeta=\displaystyle\theta\Phi_*(\xi,\tau,s)+(1\!-\!\theta)\Phi
  (\xi,\tau,s)$ for some $0<\theta<1$ (depending on $\xi,\,\tau,\,s$). Since
  there exists constant $0<c<1$ such that $\displaystyle{\Phi_*(\xi,\tau,s)
  \over\Phi(\xi,\tau,s)}\geq c$ for $\epsln$ sufficiently small, we have
  $\zeta\geq\displaystyle c\Phi(\xi,\tau,s)$. Thus
$$
\begin{array}{rcl}
  &&\displaystyle \Big|u_*'(\Phi_*(\xi,\tau,s))-u_*'(\Phi(\xi,\tau,s))\Big|
    \Big|_{\xi=\Psi(r,t,s)}\\ [0.2cm]
  &\leq &\displaystyle C\big(\Phi(\xi,\tau,s)\big)^{-\gamma}
  \big|\Phi_*(\xi,\tau,s)-\Phi(\xi,\tau,s)\big|\big|_{\xi=\Psi(r,t,s)}
    \\ [0.2cm]
  &\leq &\displaystyle C\epsln\big(e^{-\mu s}-e^{-\mu\tau}\big)
    \big(\Phi(\xi,\tau,s)\big)^{-\gamma}
    \Phi_*(\xi,\tau,s)[1-\Phi_*(\xi,\tau,s)]\big|_{\xi=\Psi(r,t,s)}
    \\ [0.2cm]
  &\leq &\displaystyle C\epsln\big(e^{-\mu s}-e^{-\mu\tau}\big)
    \big(\Phi(\xi,\tau,s)\big)^{1-\gamma}[1-\Phi(\xi,\tau,s)]
    \big|_{\xi=\Psi(r,t,s)} \\ [0.2cm]
  &=&\displaystyle C\epsln\big(e^{-\mu s}-e^{-\mu\tau}\big)
    \big(\Psi(r,t,\tau)\big)^{1-\gamma}[1-\Psi(r,t,\tau)].
\end{array}
$$
  In getting the last equality we used the following relation:
$$
  \Phi(\Psi(r,t,s),\tau,s)=\Psi(r,t,\tau) \quad
  \mbox{for}\;\; 0\leq r\leq 1, \;\; s\leq\tau\leq t.
\eqno{(5.36)}
$$
  The proof of this relation is as follows: From (5.8) we know that $\rho=
  \Psi(r,t,\tau)$ is a solution of the following problem:
$$
  \displaystyle {\partial\rho\over\partial t}+w(r,t)
  {\partial \rho\over\partial r}=0 \quad
  \mbox{for}\;\; 0\leq r\leq 1, \;\; t>\tau,\quad
  \rho|_{t=\tau}=r.
$$
  But it is easy to verify that $\rho=\Phi(\Psi(r,t,s),\tau,s)$ is also a
  solution of this problem. Hence, by uniqueness we have (5.36). Hence,
  using (5.26) and (5.14) we get
$$
\begin{array}{rcl}
  &&\displaystyle \Big|u_*'(\Phi_*(\xi,\tau,s))-{\partial w\over\partial r}
    (\Phi(\xi,\tau,s),\tau)\Big|\Big|_{\xi=\Psi(r,t,s)}\\ [0.2cm]
  &\leq &\displaystyle C\epsln\big(e^{-\mu s}-e^{-\mu\tau}\big)
    \big(\Psi_*(r,t,\tau)\big)^{1-\gamma}[1-\Psi_*(r,t,\tau)]+
    C\epsln e^{-\mu\tau}\\ [0.2cm]
  &=&\displaystyle C\epsln\big(e^{-\mu s}-e^{-\mu\tau}\big)
    \big[F_*^{-1}(F_*(r)\!+\!t\!-\!\tau)\big]^{1-\gamma}
    [1-F_*^{-1}(F_*(r)\!+\!t\!-\!\tau)]+C\epsln e^{-\mu\tau}.
\end{array}
$$
  It follows that
$$
\begin{array}{rcl}
  &&\displaystyle\int_s^t\Big|u_*'(\Phi_*(\xi,\tau,s))-{\partial w\over
  \partial r}(\Phi(\xi,\tau,s),\tau)\Big|\Big|_{\xi=\Psi(r,t,s)}d\tau
\\ [0.2cm]
  &\leq &\displaystyle C\epsln\big(e^{-\mu s}-e^{-\mu t}\big)
  \int_s^t \big[F_*^{-1}(F_*(r)\!+\!t\!-\!\tau)\big]^{1-\gamma}
  [1-F_*^{-1}(F_*(r)\!+\!t\!-\!\tau)]d\tau+
  C\epsln\int_s^t e^{-\mu\tau}d\tau\\ [0.2cm]
%  &=&\displaystyle C\epsln\big(e^{-\mu s}-e^{-\mu t}\big)
%  \int_0^{t-s}F_*^{-1}(F_*(r)\!+\!\tau')[1-F_*^{-1}(F_*(r)\!+\!\tau')]d\tau'
%  + C\epsln\big(e^{-\mu s}-e^{-\mu t}\big)\\ [0.2cm]
  &=&\displaystyle C\epsln\big(e^{-\mu s}-e^{-\mu t}\big)
  \int_{F_*(r)}^{F_*(r)+t-s}\big(F_*^{-1}(\xi)\big)^{1-\gamma}
  [1-F_*^{-1}(\xi)]d\xi+ C\epsln\big(e^{-\mu s}-e^{-\mu t}\big)
    \\ [0.2cm]
  &&\displaystyle  \qquad\qquad\qquad\qquad\qquad\qquad\qquad
    \Big(\xi=F_*(\eta),\;\;\; d\xi=F_*'(\eta)d\eta=
    {d\eta\over |u_*(\eta)|}\Big)\\ [0.2cm]
  &=&\displaystyle C\epsln\big(e^{-\mu s}-e^{-\mu t}\big)
  \int_{r}^{F_*^{-1}(F_*(r)+t-s)}{\eta^{1-\gamma}(1-\eta)\over
  |u_*(\eta)|}d\eta + C\epsln\big(e^{-\mu s}-e^{-\mu t}\big)\\ [0.2cm]
  &\leq &\displaystyle C\epsln\big(e^{-\mu s}-e^{-\mu t}\big)
   \int_0^1\eta^{-\gamma} d\eta+C\epsln\big(e^{-\mu s}-e^{-\mu t}\big)
  \\ [0.2cm]
  &=&\displaystyle C\epsln\big(e^{-\mu s}-e^{-\mu t}\big).
\end{array}
$$
  Combining this result with (5.35), we see that (5.33) follows. Finally,
  (5.34) is an immediate consequence of (5.33). $\quad\Box$
\medskip

  {\bf Corollary 5.8}\ \ {\em Under the assumption of Lemma 5.7, for $\epsln$
  sufficiently small we have
$$
  \Big|{\partial T(r,t,s)\over\partial r}-1\Big|\leq
  C\epsln\big(e^{-\mu s}-e^{-\mu t}\big), \quad
  \Big|{\partial S(\bar{r},t,s)\over\partial \bar{r}}-1\Big|\leq
  C\epsln\big(e^{-\mu s}-e^{-\mu t}\big),
$$
  and
$$
  C^{-1}\leq {\partial T(r,t,s)\over\partial r}\leq C, \quad
  C^{-1}\leq {\partial S(\bar{r},t,s)\over\partial \bar{r}}\leq C.
$$
$\Box$}

  {\bf Lemma 5.9}\ \ {\em Assume that $a\in C^1_V[0,1]$. Then we have
$$
  \|a(S(\cdot,t,s))-a\|_\infty\leq C\|a\|_1
  \epsln\big(e^{-\mu s}-e^{-\mu t}\big),
\eqno{(5.37)}
$$
  where $\|a\|_1=\max_{0\leq r\leq 1}r(1\!-\!r)|a'(r)|$. If further
  $r^2(1\!-\!r)^2a''(r)\in C[0,1]$ then we also have
$$
  \|a(S(\cdot,t,s))-a\|_{C^1_V[0,1]}\leq
  C\|a\|_2\epsln\big(e^{-\mu s}-e^{-\mu t}\big),
\eqno{(5.38)}
$$
  where  $\|a\|_2=\|a\|_1+\max_{0\leq r\leq 1}r^2(1\!-\!r)^2|a''(r)|$.}
\medskip

  {\em Proof}:\ \ We have
$$
\begin{array}{rl}
    &|a(S(r,t,s))-a(r)|=|a'(\eta)||S(r,t,s)-r|\\[0.2cm]
  \leq &\displaystyle C\epsln\eta(1\!-\!\eta)|a'(\eta)|\cdot
  {r(1-r)\over\eta(1\!-\!\eta)}\big(e^{-\mu s}-e^{-\mu t}\big)\leq
  C\|a\|_1\epsln\big(e^{-\mu s}-e^{-\mu t}\big),
\end{array}
$$
  where $\eta=(1-\theta)r+\theta S(r,t,s)$ for some $0<\theta<1$
  (depending on $r$, $t$ and $s$). In getting the last inequality we used the
  inequality
$$
  \eta(1\!-\!\eta)\geq Cr(1-r) \quad \mbox{for}\;\; 0\leq r\leq 1,
$$
  which follows from (5.25) and the following identity:
$$
  \eta(1\!-\!\eta)=(1\!-\!\theta)r(1-r)+\theta
  S(r,t,s)[1\!-\!S(r,t,s)]+
  \theta(1\!-\!\theta)[r\!-\!S(r,t,s)]^2.
$$
  Hence (5.37) is proved. Next, we compute
$$
\begin{array}{rcl}
    &&\displaystyle r(1\!-\!r)\Big|{\partial a(S(r,t,s))\over\partial r}-
    a'(r)\Big|= r(1\!-\!r)\Big|a'(S(r,t,s))
    {\partial S(r,t,s)\over\partial r}-a'(r)\Big|\\[0.3cm]
  &\leq &\displaystyle r(1\!-\!r)|a'(S(r,t,s))|\Big|
  {\partial S(r,t,s)\over\partial r}-1\Big|+
  r(1\!-\!r)|a'(S(r,t,s))-a'(r)|\\[0.3cm]
%  &\leq &\displaystyle {r(1-r)\over S(r,t,s)[1\!-\!S(r,t,s)]}\cdot
%  S(r,t,s)[1\!-\!S(r,t,s)]|a'(S(r,t,s))|\cdot
%  C\epsln\big(e^{-\mu s}-e^{-\mu t}\big)
  &\leq &\displaystyle C\|a\|_1\cdot C\epsln\big(e^{-\mu s}-e^{-\mu t}\big)+
  r(1-r)|a''(\eta)||S(r,t,s)-r|,
\end{array}
$$
  where $\eta=(1-\theta)r+\theta S(r,t,s)$ for some $0<\theta<1$
  (depending on $r$, $t$ and $s$). Similarly as before we have
$$
\begin{array}{rcl}
  r(1-r)|a''(\eta)||S(r,t,s)-r|
  &\leq & r(1-r)|a''(\eta)|\cdot
  C\epsln\big(e^{-\mu s}-e^{-\mu t}\big)r(1-r)
\\[0.3cm]
  &=&\displaystyle {r^2(1-r)^2\over\eta^2(1\!-\!\eta)^2}\cdot
  \eta^2(1\!-\!\eta)^2|a''(\eta)|\cdot C\epsln\big(e^{-\mu s}-e^{-\mu t}\big)
\\ [0.3cm]
  &\leq & C\big(\max_{0\leq r\leq 1}r^2(1\!-\!r)^2|a''(r)|\big)
  \epsln\big(e^{-\mu s}-e^{-\mu t}\big).
\end{array}
$$
  Hence (5.38) is proved. This completes the proof of Lemma 5.9  $\quad\Box$
\medskip

  {\bf Lemma 5.10}\ \ {\em Given $a\in C[0,1]$, we define a bounded linear
  operator $L$ in $C[0,1]$ by
$$
  L(q)(r)={1\over r^3}\int_0^r\!\! a(\rho)q(\rho)\rho^2 d\rho \quad
  \mbox{for}\;\; q\in C[0,1], \;\;\; 0<r\leq 1,
$$
  and $\displaystyle L(q)(0)=\lim_{r\to 0^+}L(q)(r)={1\over 3}a(0)q(0)$. Let
  $\bar{r}=T(r,t,s)$ and $r=S(\bar{r},t,s)$ be as before, and let $\widetilde{L}$
  be the following bounded linear operator in $C[0,1]$:
$$
  \widetilde{L}(q)(\bar{r})={1\over r^3}\int_0^r\!\!
  a(\rho)q(T(\rho,t,s))\rho^2 d\rho\Big|_{r=S(\bar{r},t,s)} \quad
  \mbox{for}\;\; q\in C[0,1], \;\;\; 0<\bar{r}\leq 1,
$$
  and $\displaystyle \widetilde{L}(q)(0)=\lim_{\bar{r}\to 0^+}\widetilde{L}(q)
  (\bar{r})={1\over 3}a(0)q(0)$. Assume that $a\in C^1_V[0,1]$. Then both $L$
  and $\widetilde{L}$ are bounded linear operators from $C[0,1]$ to
  $C^1_V[0,1]$, and we have
$$
  \|\widetilde{L}-L\|_{L(C[0,1],C^1_V[0,1])}\leq C\|a\|_{C^1_V[0,1]}
  \epsln\big(e^{-\mu s}-e^{-\mu t}\big).
\eqno{(5.39)}
$$
}

  {\em Proof}:\ \ We only give the proof of (5.39), because the proof of the
  assertion that both $L$ and $\widetilde{L}$ are bounded linear operators from
  $C[0,1]$ to $C^1_V[0,1]$ follows by a similar argument.

  We first note that for $q\in C[0,1]$ and $0<\bar{r}\leq 1$,
  $\widetilde{L}(q)(\bar{r})$ can be re-written as follows:
$$
  \widetilde{L}(q)(\bar{r})={1\over [S(\bar{r},t,s)]^3}\int_0^{\bar{r}}\!\!
  a(S(\rho,t,s))q(\rho)\Big[{S(\rho,t,s)\over\rho}\Big]^2
  {\partial S(\rho,t,s)\over\partial\rho}\rho^2 d\rho.
$$
  Thus
$$
\begin{array}{rcl}
  \widetilde{L}(q)(\bar{r})-L(q)(\bar{r})
  &=&\displaystyle \Big[{\bar{r}\over S(\bar{r},t,s)}\Big]^3\cdot
  {1\over \bar{r}^3}\int_0^{\bar{r}}\!\! a(S(\rho,t,s))q(\rho)
  \Big[{S(\rho,t,s)\over\rho}\Big]^2\Big[{\partial S(\rho,t,s)
  \over\partial\rho}-1\Big]\rho^2 d\rho
\\ [0.3cm]
  &&+\displaystyle \Big[{\bar{r}\over S(\bar{r},t,s)}\Big]^3\cdot
  {1\over \bar{r}^3}\int_0^{\bar{r}}\!\! a(S(\rho,t,s))q(\rho)
  \Big\{\Big[{S(\rho,t,s)\over\rho}\Big]^2-1\Big\}\rho^2 d\rho
\\ [0.3cm]
  &&+\displaystyle \Big[{\bar{r}\over S(\bar{r},t,s)}\Big]^3\cdot
  {1\over \bar{r}^3}\int_0^{\bar{r}}\!\! [a(S(\rho,t,s))-a(\rho)]
  q(\rho)\rho^2 d\rho\\ [0.3cm]
  &&+\displaystyle \Big\{\Big[{\bar{r}\over S(\bar{r},t,s)}\Big]^3-1\Big\}
  \cdot{1\over \bar{r}^3}\int_0^{\bar{r}}\!\! a(\rho)q(\rho)\rho^2 d\rho.
\end{array}
$$
  From Corollary 5.3, Lemma 5.5, Corollary 5.8 and Lemma 5.9 we know that
$$
  \Big|{\bar{r}\over S(\bar{r},t,s)}\Big|\leq C, \quad
  \Big|{\bar{r}\over S(\bar{r},t,s)}-1\Big|\leq
  C\epsln\big(e^{-\mu s}-e^{-\mu t}\big),
$$
$$
  \Big|{S(\rho,t,s)\over\rho}\Big|\leq C, \quad
  \Big|{\partial S(\rho,t,s)\over\partial\rho}-1\Big|\leq
  C\epsln\big(e^{-\mu s}-e^{-\mu t}\big),
$$
$$
  |a(S(\rho,t,s))-a(\rho)|\leq C\|a\|_1
  \epsln\big(e^{-\mu s}-e^{-\mu t}\big).
$$
  Using the above estimates, we see easily that
$$
  \max_{0\leq\bar{r}\leq 1}|\widetilde{L}(q)(\bar{r})-L(q)(\bar{r})|\leq
  C\epsln\big(e^{-\mu s}-e^{-\mu t}\big)\|a\|_{C^1_V[0,1]}\|q\|_\infty.
\eqno{(5.40)}
$$
  Next, by a simple computation we have
$$
\begin{array}{rcl}
  \bar{r}(1\!-\!\bar{r})L(q)'(\bar{r})&=&\displaystyle (1\!-\!\bar{r})
  a(\bar{r})q(\bar{r})-{3(1\!-\!\bar{r})\over \bar{r}^3}\int_0^{\bar{r}}
  a(\rho)q(\rho)\rho^2 d\rho,\\ [0.3cm]
  \bar{r}(1\!-\!\bar{r})\widetilde{L}(q)'(\bar{r})&=&\displaystyle
  {\bar{r}(1\!-\!\bar{r})\over S(\bar{r},t,s)}
  a(S(\bar{r},t,s))q(\bar{r}){\partial S(\bar{r},t,s)\over\partial \bar{r}}
\\ [0.3cm]
  &&\displaystyle -{3\bar{r}(1\!-\!\bar{r})\over [S(\bar{r},t,s)]^4}
  {\partial S(\bar{r},t,s)\over\partial \bar{r}}\int_0^{\bar{r}}
  a(S(\rho,t,s))q(\rho)\Big[{S(\rho,t,s)\over\rho}\Big]^2
  {\partial S(\rho,t,s)\over\partial\rho}\rho^2 d\rho.
\end{array}
$$
  Using these expressions and a similar argument as before we have
$$
  \sup_{0<\bar{r}<1}\bar{r}(1\!-\!\bar{r})|\widetilde{L}(q)'(\bar{r})-
  L(q)'(\bar{r})|\leq C\epsln\big(e^{-\mu s}-e^{-\mu t}\big)
  \|a\|_{C^1_V[0,1]}\|q\|_\infty.
\eqno{(5.41)}
$$
  To save spaces, we omit the details here. By (5.40) and (5.41), we see that
  (5.39) follows. $\quad\Box$
\medskip

  What we shall use later on is not (5.39), but the following immediate
  consequences of it:
$$
  \|\widetilde{L}-L\|_{L(C[0,1])}\leq C\|a\|_{C^1_V[0,1]}
  \epsln\big(e^{-\mu s}-e^{-\mu t}\big),
\eqno{(5.42)}
$$
$$
  \|\widetilde{L}-L\|_{L(C^1_V[0,1])}\leq C\|a\|_{C^1_V[0,1]}
  \epsln\big(e^{-\mu s}-e^{-\mu t}\big).
\eqno{(5.43)}
$$

\section{Decay estimates}

  In this section we establish a decay estimate for the evolution system
  $\{{\mathbb U}(t,s,V):t\geq s\geq 0\}$ obtained in Section 4, where
  $V=V(t)\in C([0,\infty),S_\epsln)$, under an additional assumption that
  $V(t)$ is exponentially decaying as $t\to\infty$.

  We first consider the special case that $V=0$. In this case we have
  ${\mathbb U}(t,s,V)=e^{(t-s){\mathbb A}(0)}$. The main result Theorem 5.1 of
  \cite{ChenCuiF} gives a decay estimate for $e^{t{\mathbb A}(0)}$ (see (6.9)
  below). But that estimate contains some singularity at $r=0$, so that it does
  not meet our requirement. In what follows we shall establish an improved
  estimate. To this end we need a preliminary lemma which gives an estimate for
  the semigroup generated by the following operator ${\mathcal L}
  ={\mathcal L}_0+a$:
$$
  {\mathcal L}q(r)=-w(r)q'(r)+a(r)q(r) \quad \mbox{for} \;\; 0\leq r\leq 1,
$$
  where $w$ and $a$ are given functions.
\medskip

  {\bf Lemma 6.1}\ \ {\em Assume that $w\in C^1[0,1]$ and satisfies $(4.1)$,
  and $a\in C^1_V[0,1]$. Then ${\mathcal L}$ generates a $C_0$ semigroup
  $e^{t{\mathcal L}}$ on $C[0,1]$ satisfying the following estimate:
$$
  \|e^{t{\mathcal L}}\|_{L(C[0,1])}\leq e^{\omega_0 t} \quad
  \mbox{for}\;\; t\geq 0,
\eqno{(6.1)}
$$
  where $\omega_0=\max_{0\leq r\leq 1}a(r)$. Moreover, $C^1_V[0,1]$ is
  ${\mathcal L}$-admissible, and for any $\omega>\omega_0$ we have
$$
  \|e^{t{\mathcal L}}\|_{L(C^1_V[0,1])}\leq C_{\omega}e^{\omega t} \quad
  \mbox{for}\;\; t\geq 0.
\eqno{(6.2)}
$$
  Here $C_\omega$ is independent of $w$ $($but depends on the constants $C_1$,
  $C_2$ in $(4.1)$ and the upper bound of $\|a\|_{C^1_V[0,1]})$.}

  {\em Proof}:\ \ By a similar argument as in the proof of Lemma 4.1 we see
  that for any $\lambda\in {\mathbb C}$ with ${\rm Re}\lambda>\omega_0$ and
  any $f\in C[0,1]$, the equation
$$
  -w(r)q'(r)+a(r)q(r)-\lambda q(r)=f(r)
$$
  has a unique solution $u\in C^1_V[0,1]$ which is given by
$$
  q(r)=e^{\int_{r_0}^r {a(\rho)-\lambda\over w(\rho)} d\rho}
  \int_r^1 {f(\eta)\over w(\eta)} e^{-\int_{r_0}^{\eta}
  {a(\rho)-\lambda\over w(\rho)} d\rho} d\eta,
$$
  where $r_0$ is an arbitrarily fixed number in $(0,1)$. Using this expression
  and a similar argument as in the proof of Lemma 4.1 we have the following
  estimate:
$$
  \max_{0\leq r\leq 1}|q(r)|\leq\max_{0\leq\eta\leq 1}{|f(\eta)|\over
  {\rm Re}\lambda-a(\eta)}\leq {1\over {\rm Re}\lambda-\omega_0}
  \max_{0\leq r\leq 1}|f(r)| \quad \mbox{for}\;\;
  {\rm Re}\lambda>\omega_0,
$$
  Hence, ${\mathcal L}$ generates a strongly continuous semigroup
  $e^{t{\mathcal L}}$ on $C[0,1]$ and the estimate (6.1) holds.

  Next, by (4.1) we see that for any $q\in C^1_V[0,1]$ we have
$$
  C_1'[\max_{0\leq r\leq 1}|q(r)|+\max_{0\leq r\leq 1}|w(r)q'(r)|]
  \leq \|q\|_{C^1_V[0,1]}\leq
  C_2'[\max_{0\leq r\leq 1}|q(r)|+\max_{0\leq r\leq 1}|w(r)q'(r)|],
\eqno{(6.3)}
$$
  where $C_1'$ and $C_2'$ are positive constants independent of $w$ (but
  depending on the constants $C_1$, $C_2$ appearing in (4.1)). Let $q_0\in
  C^1_V[0,1]$ and let $q=e^{t{\mathcal L}}q_0$. Then $q\in C([0,\infty),
  C^1_V[0,1])\cap C^1([0,\infty),C[0,1])$ and it is the solution of the
  following problem:
$$
  {\partial q\over\partial t}+w(r){\partial q\over\partial r}=a(r)q
  \quad \mbox{for}\;\; t>0, \quad q|_{t=0}=q_0.
$$
  Let $l(r,t)=w(r)\displaystyle{\partial q(r,t)\over\partial r}$. Formally
  differentiating the above equation and multiplying it with $w(r)$, we see
  that $l$ is a {\em formal solution} of the following problem:
$$
  {\partial l\over\partial t}+w(r){\partial l\over\partial r}=a(r)l+f_1(r,t),
  \quad \mbox{for}\;\; t>0, \quad l|_{t=0}=l_0,
$$
  where $f_1(r,t)=w(r)a'(r)q(r,t)$ and $l_0(r)=w(r)q_0'(r)$. Clearly, $f_1\in
  C^1([0,\infty),C[0,1])$ and $l_0\in C[0,1]$, so that by the theory of $C_0$
  semigroups (see, e. g. the discussion in Section 4.2 of \cite{Pazy};
  particularly Definition 2.3 and Theorem 2.7 there) it follows that the above
  problem has a unique so-called {\em mild solution} $l\in C([0,\infty),C[0,1]$
  and, consequently, the above formal computation makes sense, or in other words,
  $l(r,t)=w(r)\displaystyle{\partial q(r,t)\over\partial r}$ is the mild
  solution of the above problem. This means that
$$
  l(\cdot,t)=e^{t{\mathcal L}}l_0+\int_0^t e^{(t-s){\mathcal L}}f_1(\cdot,s)ds
  \quad \mbox{for}\;\; t\geq 0.
$$
  By this expression and (6.1) we have
$$
\begin{array}{rcl}
  \|l(\cdot,t)\|_\infty &\leq &\displaystyle e^{\omega_0 t}\|l_0\|_\infty+
  \int_0^t e^{\omega_0(t-s)}\|f_1(\cdot,s)\|_\infty ds \\
   &\leq &\displaystyle e^{\omega_0 t}\|l_0\|_\infty+\|wa'\|_\infty
  \int_0^t e^{\omega_0(t-s)} e^{\omega_0 s}\|q_0\|_\infty ds\\
  &\leq &\displaystyle Ce^{\omega_0 t}\|q_0\|_{C^1_V[0,1]}+
  C\|a\|_{C^1_V[0,1]}\cdot te^{\omega_0 t}\|q_0\|_\infty.
\end{array}
$$
  From this estimate and (6.1), (6.3) we immediately obtain (6.2). The proof
  is complete. $\quad\Box$
\medskip

  {\bf Lemma 6.2}\ \ {\em There exists a constant $\mu^\ast>0$ such that for
  any $0<\mu<\mu^\ast$, the semigroup $e^{t{\mathbb A}(0)}$ $(t\geq 0)$
  generated by ${\mathbb A}(0)$ satisfies the following estimate:
$$
  \|e^{t{\mathbb A}(0)}\|_{L(C[0,1])}\leq Ce^{-\mu t}\quad
  \mbox{for}\;\; t\geq 0.
\eqno{(6.1)}
$$
}

  {\em Proof}:\ \ Given $U_0=(\phi_0,\zeta_0)\in X$, let $U(t)=e^{t{\mathbb A}
  (0)}U_0=(\phi(r,t),\zeta(t))$. Then $(\phi,\zeta)$ is the unique solution of
  the following initial value problem:
$$
  \partial_t\phi+u_*(r)\partial_r\phi=a(r)\phi+{\mathcal B}(\phi)+b(r)\zeta
  \quad \mbox{for}\;\; 0\leq r\leq 1,\;\; t>0,
\eqno{(6.5)}
$$
$$
  {d\zeta\over dt}={\mathcal F}(\phi)+\kappa\zeta \quad \mbox{for}\;\; t>0,
\eqno{(6.6)}
$$
$$
  \phi(r,0)=\phi_0(r)  \quad \mbox{for}\;\; 0\leq r\leq 1, \quad
  \zeta(0)=\zeta_0.
\eqno{(6.7)}
$$
  where $a(r)$, $b(r)$, ${\mathcal B}(\phi)$, ${\mathcal F}(\phi)$ and $\kappa$
  are given in (2.24)--(2.28). By Theorem 5.1 of \cite{ChenCuiF} and the Remark
  in the end of Section 8 of \cite{ChenCuiF} we know that there exists constant
  $\sigma^\ast>0$ and a function $\hat{\phi}\in C^1(0,1]$ satisfying
$$
  \hat{\phi}(r)>0 \quad \mbox{for}\;\; 0<r\leq 1, \quad
  \hat{\phi}(r)\sim Cr^{-\theta} \quad \mbox{for}\;\; r\to 0
\eqno{(6.8)}
$$
  for some constants $1\leq\theta<3$ and $C>0$, such that the solution of the
  above problem satisfies the following estimate:
$$
  |\zeta(t)|+\sup_{0<r\leq 1}\Big|{\phi(r,t)\over\hat{\phi}(r)}\Big|
  \leq C\Big(|\zeta_0|+\sup_{0<r\leq 1}\Big|
  {\phi_0(r)\over\hat{\phi}(r)}\Big|\Big)
  (1+t)^2e^{-\sigma^\ast t} \quad \mbox{for}\;\; t\geq 0.
\eqno{(6.9)}
$$
  This particularly implies that for any $0<\sigma<\sigma^\ast$ and
  $\delta\in (0,1)$ we have
$$
  |\zeta(t)|+\sup_{\delta\leq r\leq 1}|\phi(r,t)|
  \leq C(|\zeta_0|+\sup_{0\leq r\leq 1}|\phi_0(r)|)e^{-\sigma t} \quad
  \mbox{for}\;\; t\geq 0,
\eqno{(6.10)}
$$
  because $1/\hat{\phi}(r)$ has a positive lower bound for $\delta\leq r\leq 1$
  and a finite upper bound for $0\leq r\leq 1$. In what follows we prove that
  for $\delta$ sufficiently small there also holds
$$
  \sup_{0\leq r\leq\delta}|\phi(r,t)|\leq Ce^{-\mu t} \quad
  \mbox{for}\;\;  t\geq 0
\eqno{(6.11)}
$$
  for some $\mu>0$.

  Take a nonnegative cut-off function $\varphi\in C[0,1]$ such that
$$
  \varphi(r)\leq 1 \quad \mbox{for}\;\; 0\leq r\leq 1, \quad
  \varphi(r)=1 \quad \mbox{for}\;\;  0\leq r\leq\delta, \quad
  \varphi(r)=0 \quad \mbox{for}\;\; 2\delta\leq r\leq 1.
$$
  We split ${\mathcal B}$ into a sum of two operators as follows:
$$
  {\mathcal B}(q)={\mathcal B}_1(q)+{\mathcal B}_2(q) \quad
  \mbox{for}\;\; q\in C[0,1],
\eqno{(6.12)}
$$
  where
$$
  {\mathcal B}_1(q)=-rp_*'(r)\varphi(r)\cdot{1\over r^3}
  \int_0^{\min\{r,\delta\}}g_p(\rho)q(\rho)\rho^2 d\rho,
$$
$$
\begin{array}{rl}
  {\mathcal B}_2(q)=&\displaystyle rp_*'(r)\int_0^1 g_p(\rho)
  q(\rho)\rho^2 d\rho-r^{-2}p_*'(r)[1-\varphi(r)]
  \int_0^{\min\{r,\delta\}}g_p(\rho)q(\rho)\rho^2 d\rho\\ [0.3cm]
  &-r^{-2}p_*'(r)\int_{\min\{r,\delta\}}^r g_p(\rho)q(\rho)\rho^2 d\rho,
\end{array}
$$
  and introduce
$$
  f(r,t)={\mathcal B}_2(\phi(\cdot,t))(r)+b(r)\zeta(t).
$$
  By (6.5) and the splitting (6.12), we see that $\phi$ is the solution of the
  equation
$$
  \partial_t\phi+u_*(r)\partial_r\phi=a(r)\phi+{\mathcal B}_1(\phi)
  +f(r,t) \quad \mbox{for}\;\; 0\leq r\leq 1,\;\; t>0
\eqno{(6.13)}
$$
  subject to the initial condition $\phi(r,0)=\phi_0(r)$. Introducing operators
  ${\mathcal L}(q)=-u_* q'+aq$ and ${\mathcal F}(t)=f(\cdot,t)$, we see that
  (6.13) can be rewritten as the following differential equation in $C[0,1]$:
$$
  {dq\over dt}=({\mathcal L}+{\mathcal B}_1)(q)+{\mathcal F}(t).
\eqno{(6.14)}
$$
  Using (6.8) and (6.9) we can easily show that
$$
  \|{\mathcal B}_2(\phi(\cdot,t))\|_\infty\leq C_\delta\|U_0\|
  (1+t)^2e^{-\sigma^\ast t} \quad \mbox{for}\;\; t\geq 0.
$$
  This result combined with (6.10) yields
$$
  \|{\mathcal F}(t)\|_\infty\leq C_\delta\|U_0\|e^{-\sigma t} \quad
  \mbox{for}\;\; t\geq 0.
\eqno{(6.15)}
$$
  Using the fact that $\lim_{r\to 0}rp_*'(r)=0$ and $\varphi(r)=0$ for $r\geq
  2\delta$, one can easily deduce that for any given $\epsln>0$ there exists
  corresponding $\delta>0$ such that
$$
  \|{\mathcal B}_1(q)\|_\infty\leq\epsln \|q\|_\infty.
\eqno{(6.16)}
$$
  Furthermore, from Lemma 3.1 we know that $w(r)=u_*(r)$ and $a(r)$ satisfies
  the assumptions in Lemma 6.1, so that, by Lemma 6.1, the operator ${\mathcal
  L}$ generates a strongly continuous semigroup $e^{t{\mathcal L}}$ on $C[0,1]$,
  and
$$
  \|e^{t{\mathcal L}}\|\leq e^{-\omega t}\quad \mbox{for all}\;\; t\geq 0,
\eqno{(6.17)}
$$
  where $\omega=\min_{0\leq r\leq 1}|a(r)|>0$. Since ${\mathcal B}_1$ is a bounded
  linear operator on $C[0,1]$ and, by (6.16), $\|{\mathcal B}_1\|_{L(C[0,1])}
  \leq\epsln$, it follows that the operator ${\mathcal L}+{\mathcal B}_1$ also
  generates a strongly continuous semigroup $e^{t({\mathcal L}+{\mathcal B}_1)}$
  on $C[0,1]$, and, furthermore, there holds
$$
  \|e^{t({\mathcal L}+{\mathcal B}_1)}\|\leq e^{-(\omega-\epsln) t}
  \quad \mbox{for all}\;\; t\geq 0,
\eqno{(6.18)}
$$
  In what follows we assume that $\epsln$ is sufficiently small such that
  $\omega-\epsln>0$. By (6.14) we have
$$
  q(t)=e^{t({\mathcal L}+{\mathcal B}_1)}q(0)+\int_0^t
  e^{(t-\tau)({\mathcal L}+{\mathcal B}_1)}{\mathcal F}(\tau)d\tau.
$$
  From this relation and (6.15) and (6.18) we see that for any
  $0<\mu<\min\{\sigma,\omega-\epsln\}$ there holds
$$
  \|q(t)\|_\infty\leq \|q(0)\|_\infty e^{-(\omega-\epsln) t}+
  C\|U_0\|e^{-\mu t} \quad \mbox{for}\;\; t\geq 0.
$$
  Since $q(t)=\phi(\cdot,t)$ is a solution of (6.14) with initial data $q(0)=
  \phi_0$, by this estimate we see that (6.11) follows.

  By (6.10) and (6.11), we see that (6.1) is proved. This completes the proof.
  $\quad\Box$
\medskip

  {\bf Lemma 6.3}\ \ {\em Let $\mu^\ast$ be as in Lemma 6.2. Then for any
  $0<\mu<\mu^\ast$, in addition to $(6.1)$ we also have the following estimate:
$$
  \|e^{t{\mathbb A}(0)}\|_{L(C^1_V[0,1])}\leq Ce^{-\mu t}\quad
  \mbox{for}\;\; t\geq 0.
\eqno{(6.19)}
$$
}

  {\em Proof}:\ \ We first show that ${\mathcal B}$ and ${\mathcal F}$ satisfy
  the following properties: For any $q\in C^1_V[0,1]$,
$$
  \|{\mathcal B}(u_*q')\|_\infty+\|u_*{\mathcal B}(q)'\|_\infty
  \leq C\|q\|_\infty,
\eqno{(6.20)}
$$
$$
  |{\mathcal F}(u_*q')|\leq C\|q\|_\infty.
\eqno{(6.21)}
$$

  Using the facts that $\lim_{r\to 0}rp_*'(r)=0$ and $\lim_{r\to 0}r^2p_*''(r)
  =0$ (see (3.3)) we can easily prove that
$$
  \|u_*{\mathcal B}(q)'\|_\infty\leq C\|q\|_\infty \quad \mbox{for}
  \;\; q\in C[0,1].
$$
  To estimate $\|{\mathcal B}(u_*q')\|_\infty$ we compute:
$$
  {1\over r^3}\int_0^rg_p(\rho)u_*(\rho)q'(\rho)\rho^2 d\rho
  ={1\over r}u_*(r)g_p(r)q(r)-{1\over r^3}\int_0^r m(\rho)q(\rho)\rho^2 d\rho.
$$
  where $m(\rho)=\displaystyle g_p'(\rho)u_*(\rho)+g_p(\rho)u_*'(\rho)+
  {2\over\rho}u_*(\rho)g_p(\rho)$. Taking $r=1$ we particularly obtain
$$
  \int_0^1 g_p(\rho)u_*(\rho)q'(\rho)\rho^2 d\rho
  =u_*(1)g_p(1)q(1)-\int_0^1 m(\rho)q(\rho)\rho^2 d\rho.
$$
  Since $g_p\in C^1[0,1]$, $u_*\in C^1[0,1]$ and $u_*(0)=0$, we see that
  $\displaystyle {1\over r}u_* g_p$ and $m$ both belong to $C[0,1]$. Hence,
  from the above expressions we see immediately that
$$
  \|{\mathcal B}(u_*q')\|_\infty=\sup_{0<r<1}\Big|rp_*'(r)
  \Big[\int_0^1 g_p(\rho)u_*(\rho)q'(\rho)\rho^2 d\rho-
  {1\over r^3}\int_0^r g_p(\rho)u_*(\rho)q'(\rho)\rho^2 d\rho\Big]\Big|
  \leq C\|q\|_\infty.
$$
  Similarly we also have
$$
  |{\mathcal F}(u_*q')|=\Big|\int_0^1 g_p(\rho)u_*(\rho)
  q'(\rho)\rho^2 d\rho\Big|\leq C\|q\|_\infty.
$$
  This verifies (6.20) and (6.21).

  We now proceed to prove (6.19). Let $U_0\in X_0$ and $U=e^{t{\mathbb A}(0)}
  U_0$. From the proof of Lemma 4.2 we know that $U\in C([0,\infty),X_0)\cap
  C^1([0,\infty),X)$. Let $U_0=(q_0,y_0)$ and $U=(q,y)$. Then $(q,y)$ is the
  solution of the following problem:
$$
\left\{
\begin{array}{l}
  \displaystyle {\partial q\over\partial t}+u_*(r){\partial q\over\partial r}
  =a(r)q+{\mathcal B}(q)+b(r)y
  \quad \mbox{for} \;\; 0\leq r\leq 1, \;\; t>0,\\ [0.1cm]
  \displaystyle {dy\over dt}={\mathcal F}(q)+\kappa y
  \quad \mbox{for} \;\; t>0,\\ [0.1cm]
  q|_{t=0}=q_0(r)\quad \mbox{for} \;\; 0\leq r\leq 1, \quad
  \mbox{and}\quad y|_{t=0}=y_0.
\end{array}
\right.
$$
  Let $l(r,t)=u_*(r)\displaystyle{\partial q(r,t)\over\partial r}$. As in the
  proof of Lemma 6.1, by formally differentiating the first equation above in
  $r$ and multiplying it with $u_*(r)$, we see that $(l,y)$ is a ``formal
  solution'' of the following problem:
$$
\left\{
\begin{array}{l}
  \displaystyle {\partial l\over\partial t}+u_*(r){\partial l\over\partial r}
  =a(r)l+{\mathcal B}(l)+b(r)y+f_1(r,t)
  \quad \mbox{for} \;\; 0\leq r\leq 1, \;\; t>0,\\ [0.1cm]
  \displaystyle {dy\over dt}={\mathcal F}(l)+\kappa y+c_1(t)
  \quad \mbox{for} \;\; t>0,\\ [0.1cm]
  l|_{t=0}=l_0(r)\quad \mbox{for} \;\; 0\leq r\leq 1, \quad
  \mbox{and}\quad y|_{t=0}=y_0,
\end{array}
\right.
$$
  where $l_0(r)=u_*(r)q_0'(r)$, $c_1(t)=\displaystyle{\mathcal F}(q)-
  {\mathcal F}\big(u_*{\partial q\over\partial r}\big)$, and
$$
  f_1(r,t)=\displaystyle u_*(r)a'(r)q(r,t)-
  {\mathcal B}\Big(u_*{\partial q\over\partial r}\Big)
  +u_*(r){\partial{\mathcal B}(q)\over\partial r}+[u_*(r)b'(r)-b(r)]y(t).
$$
  We denote $W(t)=(l(\cdot,t),y(t))$, $W_0=(l_0,y_0)$ and $F_1(t)=
  (f_1(\cdot,t),c_1(t))$. Then the above problem can be rewritten as follows:
$$
  {dW\over dt}={\mathbb A}(0)W+F_1(t) \quad \mbox{for}\;\; t>0, \quad
  W(0)=W_0.
$$
  Using the fact that $U\in C^1([0,\infty),X)$, Corollary 3.2 and (6.20),
  (6.21), we can easily prove that $F_1\in C^1([0,\infty),X)$. Thus, by a
  similar argument as in the proof of Lemma 6.1 we see that the above formal
  computation makes sense and $W=(l,y)=(u_*\displaystyle{\partial q\over
  \partial r},y)$ is the unique mild solution of the above problem, which
  means that
$$
  W(t)=e^{t{\mathbb A}(0)}W_0+\int_0^t e^{(t-s){\mathbb A}(0)}F_1(s)ds
  \quad \mbox{for}\;\; t\geq 0.
$$
  It follows by Lemma 6.2 that for any given $0<\mu<\mu_*$ we have
$$
  \|W(t)\|_X\leq Ce^{-\mu t}\|W_0\|_X+C\int_0^t e^{-\mu(t-s)}\|F_1(s)\|_X ds
  \quad \mbox{for}\;\; t\geq 0.
$$
  Using (6.20), (6.21) and the fact that $\|U(t)\|_X\leq Ce^{-\mu t}\|U_0\|_X$
  ensured by (6.1), we see that
$$
  \|F_1(t)\|_X\leq C\|U(t)\|_X\leq Ce^{-\mu t}\|U_0\|_X
    \quad \mbox{for}\;\; t\geq 0.
$$
  Hence, by a similar argument as in the proof of Lemma 6.1 we obtain
$$
  \|U(t)\|_{X_0}\leq C(1+t)e^{-\mu t}\|U_0\|_{X_0}
    \quad \mbox{for}\;\; t\geq 0.
$$
  Now, for any given $0<\mu<\mu_*$ we arbitrarily take a $\bar{\mu}\in(\mu,
  \mu_*)$ and first use the above estimate to $\bar{\mu}$ and next use the
  elementary inequality $(1+t)e^{-\bar{\mu}t}\leq Ce^{-\mu t}$, we see that
  (6.19) follows. This completes the proof. $\quad\Box$
\medskip

  In the sequel we consider the evolution system ${\mathbb U}(t,s, V)$ for a
  general $V=V(t)\in C([0,\infty),X)$ satisfying the following condition:
  For some positive constants $\bar{\mu}$, $\epsln$ and $C_0$,
$$
  \|V(t)\|_{X}\leq C_0\epsln e^{-\bar{\mu}t} \quad \mbox{for} \;\; t\geq 0.
\eqno{(6.22)}
$$

  {\bf Lemma 6.4}\ \ {\em Assume that $V=V(t)\in C([0,\infty),X)$ and it
  satisfies $(6.22)$. Let $\mu^*$ be as in Lemma 6.2. Then for any $0<\mu<
  \mu^*$ there exists corresponding $\epsln_0>0$ $($depending on $\mu$,
  $\bar\mu$ and $C_0)$ such that if $0<\epsln\leq\epsln_0$ then the following
  estimates hold:
$$
  \|{\mathbb U}(t,s, V)\|_{L(X)}\leq C_1 e^{-\mu t} \quad
  \mbox{for} \;\; t\geq 0,
\eqno{(6.23)}
$$
$$
  \|{\mathbb U}(t,s, V)\|_{L(X_0)}\leq C_2 e^{-\mu t} \quad
  \mbox{for} \;\; t\geq 0,
\eqno{(6.24)}
$$
  where $C_1$ and $C_2$ are positive constants depending only on $\mu$ and
  independent of $\bar\mu$ and $C_0$.}
\medskip

  {\em Proof}:\ \ Given $0<\mu<\mu^*$ we take a $\mu_1\in (\mu,\mu^*)$ and fix
  it. By Lemmas 6.2 and 6.3, we have the following estimates:
$$
  \|e^{t{\mathbb A}(0)}\|_{L(C[0,1])}\leq C_1 e^{-\mu_1 t}\quad
  \mbox{for}\;\; t\geq 0,
\eqno{(6.25)}
$$
$$
  \|e^{t{\mathbb A}(0)}\|_{L(C^1_V[0,1])}\leq C_2 e^{-\mu_1 t}\quad
  \mbox{for}\;\; t\geq 0.
\eqno{(6.26)}
$$
  Let $U_0=(q_0,s_0)$ be an arbitrary point in $X$, and let $U={\mathbb U}(t,s,V)
  U_0$. By definition, $U$ is the solution of the problem (4.13). Let $U=(q,y)$.
  Then (4.13) can be rewritten as follows:
$$
\left\{
\begin{array}{l}
  \displaystyle {\partial q\over\partial t}+w_V(r,t){\partial q\over\partial r}
  =a(r)q+{\mathcal B}q+b(r)y
  \quad \mbox{for} \;\; 0\leq r\leq 1, \;\; t>s,\\ [0.1cm]
  \displaystyle {dy\over dt}={\mathcal F}(q)+\kappa y
  \quad \mbox{for} \;\; t>s,\\ [0.1cm]
  q|_{t=s}=q_0(r)\quad \mbox{for} \;\; 0\leq r\leq 1, \quad
  \mbox{and}\quad y|_{t=s}=y_0,
\end{array}
\right.
\eqno{(6.27)}
$$
  Let $\widetilde{q}(\bar{r},t,s)=q(S(\bar{r},t,s),t)$ or $q(r,t)=\widetilde{q}
  (T(r,t,s),t,s)$. Then by using (5.9) we see that (6.27) is transformed into
  the following problem:
$$
\left\{
\begin{array}{l}
  \displaystyle {\partial\widetilde{q}\over\partial t}+
  u_*(\bar{r}){\partial\widetilde{q}\over\partial \bar{r}}
  =\widetilde{a}(\bar{r},t,s)\widetilde{q}+
  \widetilde{\mathcal B}\widetilde{q}+\widetilde{b}(\bar{r},t,s)s
  \quad \mbox{for} \;\; 0\leq \bar{r}\leq 1, \;\; t>s,\\ [0.1cm]
  \displaystyle {ds\over dt}=\widetilde{\mathcal F}(\widetilde{q})(t,s)+
  \kappa s \quad \mbox{for} \;\; t>s,\\ [0.1cm]
  \widetilde{q}|_{t=s}=q_0(\bar{r})\quad \mbox{for} \;\; 0\leq \bar{r}\leq 1, \quad
  \mbox{and}\quad s|_{t=s}=s_0,
\end{array}
\right.
\eqno{(6.28)}
$$
  where $\widetilde{a}(\bar{r},t,s)=a(S(\bar{r},t,s))$, $\widetilde{b}
  (\bar{r},t,s)=b(S(\bar{r},t,s))$,
$$
  \widetilde{\mathcal B}\widetilde{q}=
  \displaystyle rp_*'(r)\Big[\int_0^1 g_p(\rho)
  \widetilde{q}(T(\rho,t,s),t,s)\rho^2 d\rho-{1\over r^3}
  \int_0^r g_p(\rho)\widetilde{q}(T(\rho,t,s),t,s)\rho^2 d\rho\Big]
  \Big|_{r=S(\bar{r},t,s)},
$$
  and $\widetilde{\mathcal F}(\widetilde{q})(t,s)\displaystyle
  =\int_0^1 g_p(\rho)\widetilde{q}(T(\rho,t,s),t,s)\rho^2 d\rho$.
  We define a family of bounded linear operators $\widetilde{{\mathbb B}}
  (t,s,V):X\to X$ ($t\geq s\geq 0$) as follows:
$$
  \widetilde{{\mathbb B}}(t,s,V)=\left(
\begin{array}{cc}
  \widetilde{a}(\cdot,t,s)+\widetilde{\mathcal B}
  & \quad \widetilde{b}(\cdot,t,s)\\
  \widetilde{\mathcal F} & \quad \kappa
\end{array}
\right).
$$
  We also denote $\widetilde{U}=(\widetilde{q},y)$. Then (6.28) can be
  rewritten as follows:
$$
\left\{
\begin{array}{l}
  \displaystyle {d\widetilde{U}\over dt}=
  {\mathbb A}_0(U_*)\widetilde{U}+\widetilde{{\mathbb B}}(t,s,V)\widetilde{U}
  \quad \mbox{for} \;\; t>s,\\ [0.1cm]
  \widetilde{U}|_{t=s}=U_0.
\end{array}
\right.
\eqno{(6.29)}
$$
  Recalling that ${\mathbb A}(0)={\mathbb A}_0(U_*)+{\mathbb B}$ and denoting
$$
  \widetilde{{\mathbb E}}(t,s,V)=\widetilde{{\mathbb B}}(t,s,V)-{\mathbb B}
  =\left(
\begin{array}{cc}
  \widetilde{a}(\cdot,t,s)-a+\widetilde{\mathcal B}-{\mathcal B}
  & \quad \widetilde{b}(\cdot,t,s)-b\\
  \widetilde{\mathcal F}-{\mathcal F} & \quad 0
\end{array}
\right),
$$
  we see that
$$
  {\mathbb A}_0(U_*)+\widetilde{{\mathbb B}}(t,s,V)=
  {\mathbb A}_0(U_*)+{\mathbb B}+\widetilde{{\mathbb E}}(t,s,V)=
  {\mathbb A}(0)+\widetilde{{\mathbb E}}(t,s,V).
$$
  Hence, (6.29) can be further rewritten as follows:
$$
\left\{
\begin{array}{l}
  \displaystyle {d\widetilde{U}\over dt}=
  {\mathbb A}(0)\widetilde{U}+\widetilde{{\mathbb E}}(t,s,V)\widetilde{U}
  \quad \mbox{for} \;\; t>s,\\ [0.1cm]
  \widetilde{U}|_{t=s}=U_0.
\end{array}
\right.
\eqno{(6.30)}
$$
  We know that (6.30) is equivalent to the following integral equation:
$$
  \widetilde{U}(t,s)=e^{(t-s){\mathbb A}(0)}U_0+
  \int_s^t e^{(t-\tau){\mathbb A}(0)}\widetilde{{\mathbb E}}(\tau,s,V)
  \widetilde{U}(\tau,s) \quad \mbox{for} \;\; t\geq s.
\eqno{(6.31)}
$$
  By Corollary 3.2 and Lemma 5.9 we have
$$
  \|\widetilde{a}(\cdot,t,s)-a\|_\infty\leq
  \|\widetilde{a}(\cdot,t,s)-a\|_{C^1_V[0,1]}\leq
  C\epsln,
$$
$$
  \|\widetilde{b}(\cdot,t,s)-b\|_\infty\leq
  \|\widetilde{b}(\cdot,t,s)-b\|_{C^1_V[0,1]}\leq
  C\epsln,
$$
  and by Corollary 3.2, Lemma 5.9 and Lemma 5.10 we have
$$
  \|\widetilde{\mathcal B}-{\mathcal B}\|_{L(C[0,1])}\leq
  C\epsln,\quad
  \|\widetilde{\mathcal B}-{\mathcal B}\|_{L(C^1_V[0,1])}\leq
  C\epsln,
$$
$$
  \|\widetilde{\mathcal F}-{\mathcal F}\|_{L(C[0,1],{\mathbb R})}
  \leq C\epsln.
$$
  It follows that
$$
  \|\widetilde{{\mathbb E}}(\tau,s,V)\|_{L(X)}\leq C\epsln,
  \quad
  \|\widetilde{{\mathbb E}}(\tau,s,V)\|_{L(X_0)}\leq C\epsln.
\eqno{(6.32)}
$$
  From (6.25), (6.26), (6.31) and (6.32) we obtain:
$$
  \|\widetilde{U}(t,s)\|_{X}\leq C_1 e^{-\mu_1(t-s)}\|U_0\|_{X}+
  C\epsln\int_s^t e^{-\mu_1(t-\tau)}\|\widetilde{U}(\tau,s)\|_{X},
$$
$$
  \|\widetilde{U}(t,s)\|_{X_0}\leq C_2 e^{-\mu_1(t-s)}\|U_0\|_{X_0}+
  C\epsln\int_s^t e^{-\mu_1(t-\tau)}\|\widetilde{U}(\tau,s)\|_{X_0}.
$$
  By Gronwall lemma, these inequalities yield
$$
  \|\widetilde{U}(t,s)\|_{X}\leq C_1e^{-(\mu_1-C\epsln)t}\|U_0\|_{X},
$$
$$
  \|\widetilde{U}(t,s)\|_{X_0}\leq C_2e^{-(\mu_1-C\epsln)t}\|U_0\|_{X_0}.
$$
  Hence, by taking $\epsln$ sufficiently small such that $\mu_1-C\epsln\geq
  \mu$, we obtain (6.23) and (6.24). This completes the proof. $\quad\Box$
\medskip

\section{The proof of Theorem 1.1}

  In order to prove Theorem 1.1, we let $\mu^*$ be as in Lemma 6.2 and
  arbitrarily fix a number $0<\mu<\mu^*$. Let $\epsln$ be a positive number to
  be specified later. For any fixed $U_0\in X_0$ satisfying $\|U_0\|_{X_0}\leq
  \epsln$, we denote by ${\mathbf M}$ the set of all functions $V=V(t)\in
  C([0,\infty),X)$ satisfying the following conditions:
$$
  V(0)=U_0,\quad  \|V(t)\|_{X}\leq 2C_1\epsln e^{-\mu t}\;\;
  \mbox{for}\;\; t\geq 0,
\eqno{(7.1)}
$$
  where $C_1$ is the constant appearing in (6.23). We introduce a metric $d$ on
  ${\mathbf M}$ by defining
$$
  d(V_1,V_2)=\sup_{t\geq 0} e^{\mu t}\|V_1(t)-V_2(t)\|_X \quad
  \mbox{for} \;\; V_1,V_2\in {\mathbf M}.
$$
  It is evident that $({\mathbf M},d)$ is a complete metric space. Given $V\in
  {\mathbf M}$, we consider the following initial value problem:
$$
\left\{
\begin{array}{l}
  \displaystyle{dU(t)\over dt}={\mathbb A}(V(t))U(t)+{\mathbb G}(U(t))
  \quad \mbox{for}\;\; t>0,\\
  U(0)=U_0.
\end{array}
\right.
\eqno{(7.2)}
$$
\medskip

  {\bf Lemma 7.1}\ \ {\em If $\epsln$ is sufficiently small then for any $V\in
  {\mathbf M}$ the problem $(7.2)$ has a unique solution $U\in C([0,\infty),
  X_0)\cap C^1([0,\infty),X)$ which satisfies the following estimates:
$$
  \|U(t)\|_{X}\leq 2C_1\epsln e^{-\mu t},\quad
  \|U(t)\|_{X_0}\leq C\epsln e^{-\mu t}, \quad
  \|U'(t)\|_X\leq C\epsln e^{-\mu t} \quad \mbox{for}\;\; t\geq 0,
\eqno{(7.3)}
$$
  where $C_1$ is as before, and $C$ is another constant independent of $V$.}
\medskip

  {\em Proof}:\ \ We denote
$$
  \widetilde{\mathbf M}=\{U\in C([0,\infty),X_0):
  \|U(t)\|_{X}\leq 2C_1\epsln e^{-\mu t}\;\;\mbox{and} \;\;
  \|U(t)\|_{X_0}\leq 2C_2\epsln e^{-\mu t} \;\; \mbox{for}\;\; t\geq 0\},
$$
  and introduce a metric $d$ on it by defining
$$
  d(U_1,U_2)=\sup_{t\geq 0}e^{\mu t}\|U_1(t)-U_2(t)\|_{X_0}\quad
  \mbox{for} \;\; U_1,U_2\in \widetilde{\mathbf M}.
$$
  Here $C_1$ and $C_2$ are positive constants appearing in (6.23) and (6.24),
  respectively. $(\widetilde{\mathbf M},d)$ is clearly a complete metric space.
  Given $U\in \widetilde{\mathbf M}$, we consider the following initial value
  problem:
$$
\left\{
\begin{array}{l}
  \displaystyle{d\widetilde{U}(t)\over dt}={\mathbb A}(V(t))\widetilde{U}(t)+
  {\mathbb G}(U(t)) \quad \mbox{for}\;\; t>0,\\
  \widetilde{U}(0)=U_0.
\end{array}
\right.
\eqno{(7.4)}
$$
  Since $U(t)\in C([0,\infty),X_0)$, by Corollary 3.3 we have ${\mathbb G}(U(t))
  \in C([0,\infty),X_0)$. It follows by Corollary 4.4 that the above problem
  has a unique solution $\widetilde{U}\in C([0,\infty),X_0)\cap
  C^1([0,\infty),X)$, and is given by
$$
  \widetilde{U}(t)={\mathbb U}(t,0,V)U_0+
  \int_0^t{\mathbb U}(t,s,V){\mathbb G}(U(s))ds.
\eqno{(7.5)}
$$
  Using this expression and Lemma 6.4 and (2.23) we have
$$
\begin{array}{rcl}
  \|\widetilde{U}(t)\|_{X}&\leq &\displaystyle C_1 e^{-\mu t}\|U_0\|_{X}+
  C_1\int_0^t e^{-\mu(t-s)}\|{\mathbb G}(U(s))\|_{X}ds\\
  &\leq &\displaystyle C_1\epsln e^{-\mu t}+
  C\int_0^t e^{-\mu(t-s)}\|U(s)\|_{X}^2ds\\
  &\leq &\displaystyle C_1\epsln e^{-\mu t}+
  C\epsln^2\int_0^t e^{-\mu(t-s)}e^{-2\mu s}ds\\
  &\leq &\displaystyle C_1\epsln e^{-\mu t}+
  C\epsln^2 e^{-\mu t}\leq 2C_1\epsln e^{-\mu t}.
\end{array}
$$
  The last inequality holds when $\epsln$ is sufficiently small. Similarly, by
  using Lemma 6.4 and (3.9) we also have
$$
  \|\widetilde{U}(t)\|_{X_0}\leq 2C_2\epsln e^{-\mu t},
$$
  when $\epsln$ is sufficiently small. Hence $\widetilde{U}\in
  \widetilde{\mathbf M}$. We now define a mapping $\widetilde{\mathbf
  S}:\widetilde{\mathbf M}\to\widetilde{\mathbf M}$ by setting
  $\widetilde{\mathbf S}(U)=\widetilde{U}$ for every $U\in \widetilde{\mathbf
  M}$. We claim that $\widetilde{\mathbf S}$ is a contraction mapping.
  Indeed, for any $U_1,U_2\in\widetilde{\mathbf M}$ let
  $\widetilde{U}_1={\mathbf S}(U_1)$, $\widetilde{U}_2={\mathbf S}(U_2)$ and
  $W=\widetilde{U}_1-\widetilde{U}_2$. Then $W$ satisfies
$$
\left\{
\begin{array}{l}
  \displaystyle{dW(t)\over dt}={\mathbb A}(V(t))W(t)+
  [{\mathbb G}(U_1(t))-{\mathbb G}(U_2(t))] \quad \mbox{for}\;\; t>0,\\
  W(0)=0,
\end{array}
\right.
$$
  so that
$$
  W(t)=\int_0^t {\mathbb U}(t,s,V)
  [{\mathbb G}(U_1(s))-{\mathbb G}(U_2(s))]ds.
$$
  It follows by a similar argument as before that
$$
\begin{array}{rcl}
  \|W(t)\|_{X_0}&\leq &\displaystyle C_2\int_0^t e^{-\mu(t-s)}
  \|{\mathbb G}(U_1(s)-{\mathbb G}(U_2(s)))\|_{X_0}ds\\
  &\leq &\displaystyle C_2\int_0^t e^{-\mu(t-s)}\|U_1(s)-U_2(s)\|_{X_0}ds
  \Big(\int_0^1\|{\mathbb G}'(\theta U_1(s)+(1-\theta)U_2(s))\|_{L(X_0)}
  d\theta\Big)\\
  &\leq &\displaystyle C\int_0^t e^{-\mu(t-s)}\|U_1(s)-U_2(s)\|_{X_0}ds
  \Big(\int_0^1\|\theta U_1(s)+(1-\theta)U_2(s)\|_{X_0} d\theta\Big)\\
  &\leq &\displaystyle C\int_0^t e^{-\mu(t-s)}\cdot d(U_1,U_2)
  e^{-\mu s}\cdot 2C_2\epsln e^{-\mu s} ds\leq C\epsln d(U_1,U_2) e^{-\mu t}.
\end{array}
$$
  Thus for $\epsln$ sufficiently small we have
$$
   d(\widetilde{U}_1,\widetilde{U}_2)=
   \sup_{t\geq 0}e^{\mu t}\|W(t)\|_{X_0}\leq
   {1\over 2} d(U_1,U_2),
$$
  showing that $\widetilde{\mathbf S}$ is a contraction mapping, as we claimed.
  Thus, by the Banach fixed point theorem we see that $\widetilde{\mathbf S}$
  has a unique fixed point in $\widetilde{\mathbf M}$, which is clearly a
  solution of the problem (7.2) in $C([0,\infty),X_0)$. Uniqueness of the
  solution follows from a standard argument.

  From the above argument we see that the solution $U$ of (7.2) satisfies the
  first two inequalities in (7.3), and $U\in C^1([0,\infty),X)$. It remains to
  prove that $U$ also satisfies the last inequality in (7.3). The argument is
  as follows. First, it is straightforward to deduce from the condition (7.1)
  that for sufficiently small $\epsln>0$, we have $w_V(r,t)/r(1-r)\in C[0,1]$
  and there exist positive constants $C_1$ and $C_2$ independent of $V$ such
  that
$$
  -C_1 r(1-r)\leq w_V(r,t)\leq -C_2 r(1-r) \quad
  \mbox{for}\;\; 0\leq r\leq 1\;\; \mbox{and} \;\; t\geq 0.
\eqno{(7.6)}
$$
  It follows that for the solution $U=(q,s)$ of (7.2) we have
$$
  \sup_{0\leq r\leq 1}\Big|w_V(r,t){\partial q(r,t)\over\partial r}\Big|
  \leq C\sup_{0\leq r\leq 1}\Big|r(1-r){\partial q(r,t)\over\partial r}\Big|
  \leq C\|q(\cdot,t)\|_{C^1_V[0,1]}.
$$
  Using this result and the equation (7.2) we see that
$$
  \|U'(t)\|_X\leq \|{\mathbb A}(V(t))U(t)\|_X+\|{\mathbb G}(U(t)\|_X
  \leq C\|U(t)\|_{X_0}+C\|U(t)\|_X^2\leq C\epsln e^{-\mu t}
$$
  for all $t\geq 0$. This completes the proof of Lemma 7.1. $\quad\Box$
\medskip

  Lemma 7.1 particularly implies that for every $V$ in ${\mathbf M}$, the
  solution $U$ of (7.2) also belongs to ${\mathbf M}$. Thus we can define a
  mapping ${\mathbf S}:{\mathbf M}\to {\mathbf M}$ as follows: For any $V\in
  {\mathbf M}$,
$$
  {\mathbf S}(V)=U=\mbox{the solution of (7.2)}.
$$
%\medskip

  {\bf Lemma 7.2}\ \ {\em For $\epsln$ sufficiently small, ${\mathbf S}$ is a
  contraction mapping.}
\medskip

  {\em Proof}:\ \ Let $V_1,V_2\in {\mathbf M}$ and denote $U_1={\mathbf S}(V_1)$,
  $U_2={\mathbf S}(V_2)$ and $W=U_1-U_2$. Then $W$ satisfies:
$$
\left\{
\begin{array}{l}
  \displaystyle{dW(t)\over dt}={\mathbb A}(V_1(t))W(t)+[{\mathbb A}(V_1(t))-
  {\mathbb A}(V_2(t))]U_2(t)+[{\mathbb G}(U_1(t))-{\mathbb G}(U_2(t))]
  \quad \mbox{for}\;\; t>0,\\
  W(0)=0.
\end{array}
\right.
$$
  Thus
$$
  W(t)=\int_0^t {\mathbb U}(t,s,V_1)[{\mathbb A}(V_1(s))-
  {\mathbb A}(V_2(s))]U_2(s)ds+\int_0^t {\mathbb U}(t,s,V_1)
  [{\mathbb G}(U_1(s))-{\mathbb G}(U_2(s))]ds.
\eqno{(7.7)}
$$
  Since the first component of $[{\mathbb A}(V_1(s))-{\mathbb A}(V_2(s))]
  U_2(s)$ is equal to $[w_{V_2}(r,s)-w_{V_1}(r,s)]q_2'(r,s)$ and the second
  component is zero, we have
$$
\begin{array}{rcl}
  &&\|[{\mathbb A}(V_1(s))-{\mathbb A}(V_2(s))]U_2(s)\|_X
  =\displaystyle\max_{0\leq r\leq 1}|[w_{V_1}(r,s)-
  w_{V_2}(r,s)]q_2'(r,s)|\\
  &\leq &\displaystyle\sup_{0\leq r\leq 1}\Big|{w_{V_1}(r,s)-
  w_{V_2}(r,s)\over r(1-r)}\Big|\max_{0\leq r\leq 1}|r(1\!-\!r)q_2'(r,s)|
  \leq \displaystyle C\|V_1(s)-V_2(s)\|_X\|U_2(s)\|_{X_0}.
\end{array}
$$
  Besides, from (2.23) we have
$$
\begin{array}{rcl}
  \|{\mathbb G}(U_1(s))-{\mathbb G}(U_2(s))\|_X &=&\displaystyle
  \Big\|\int_0^1 {\mathbb G}'(\theta U_1(s)+(1\!-\!\theta)U_2(s))
  [U_1(s)-U_2(s)]d\theta\Big\|_X\\
  &\leq &\displaystyle C\big(\|U_1(s)\|_X+\|U_2(s)\|_X\big)
  \|U_1(s)-U_2(s)\|_X.
\end{array}
$$
  Thus, by (7.7) and Lemma 6.4 we have
$$
\begin{array}{rcl}
  \|U_1(t)-U_2(t)\|_X &\leq &\displaystyle C_1\int_0^t e^{-\mu (t-s)}
   \|[{\mathbb A}(V_1(s))-{\mathbb A}(V_2(s))]U_2(s)\|_X ds\\
  &&\displaystyle +C_1\int_0^t e^{-\mu (t-s)}
  \|[{\mathbb G}(U_1(s))-{\mathbb G}(U_2(s))]\|_X ds\\
  &\leq &\displaystyle C\int_0^t e^{-\mu (t-s)}\|V_1(s)-V_2(s)\|_X
  \|U_2(s)\|_{X_0}ds\\
  &&\displaystyle +C\int_0^t e^{-\mu (t-s)}\big(\|U_1(s)\|_X+\|U_2(s)\|_X\big)
  \|U_1(s)-U_2(s)\|_X ds\\
  &\leq &\displaystyle C\sup_{s\geq 0} e^{\mu s}\|V_1(s)-V_2(s)\|_X\cdot
  \sup_{s\geq 0} e^{\mu s}\|U_2(s)\|_{X_0}\int_0^t e^{-\mu (t-s)}\cdot
  e^{-2\mu s}ds\\
  &&\displaystyle +C\sup_{s\geq 0} e^{\mu s}\big(\|U_1(s)\|_X+\|U_2(s)\|_X\big)
  \\
  &&\displaystyle \quad\cdot\sup_{s\geq 0} e^{\mu s}
  \|U_1(s)-U_2(s)\|_{X}\int_0^t e^{-\mu (t-s)}\cdot e^{-2\mu s}ds\\
  &\leq &\displaystyle C\epsln e^{-\mu t}d(V_1,V_2)+
   C\epsln e^{-\mu t}d(U_1,U_2).
\end{array}
$$
  Therefore,
$$
  d(U_1,U_2)=\sup_{t\geq 0} e^{\mu t}\|U_1(t)-U_2(t)\|_X\leq
  C\epsln d(V_1,V_2)+ C\epsln d(U_1,U_2),
$$
  by which the desired assertion immediately follows. This completes the proof
  of Lemma 7.2. $\quad\Box$
\medskip

  By Lemma 7.2, if $\epsln$ is sufficiently small then the mapping ${\mathbf S}$
  has a unique fixed point $U$ in ${\mathbf M}$. Clearly, $U$ is a global
  solution of the equation (2.22) subject to the initial condition $U(0)=U_0$.
  Moreover, by Lemma 7.1 we know that the image of ${\mathbf S}$ is contained
  in $\widetilde{\mathbf M}$, so that $U$ satisfies (7.3). From this result all
  assertions of Theorem 1.1 easily follows. The proof of Theorem 1.1 is
  complete.
\medskip

  {\bf Acknowledgement}\ \ This work is supported by the China National Natural
  Science Foundation under the Grant number 10471157. Part of this work was
  prepared when the author was visiting the Ecole Normale Superieur (ENS)
  during June 1 $\sim$ August 28, 2007. He wishes to acknowledge his sincere
  thanks to the Department of Mathematics and Applications of ENS for
  hospitality and the French Ministry of Foreign Affairs, particularly the
  Section of Science and Technology of the French Consulate at Guangzhou, for
  financial support to his visit.

\end{document}